\documentclass[aap,MSNbibl,citesort,seceqn,dvips]{arximspdf}
\usepackage{mathrsfs}
\usepackage{graphicx}

\doi{10.1214/10-AAP700}
\volume{21}
\issue{2}
\pubyear{2011}
\firstpage{546}
\lastpage{588}

\makeatletter
\newcommand{\Lvb}{L\'{e}vy }

\newcommand{\cadlag}{c\`{a}dl\`{a}g }
\newcommand{\Ito}{It\^{o} }

\newcommand{\var}{ \mbox{Var} }

\newcommand{\cov}{ \operatorname{Cov} }

\newcommand{\refeq}[1]{(\ref{#1})}

\newtheorem{thms}{Theorem}[section]
\newtheorem{cor}{Corollary}[section]
\newtheorem{lema}{Lemma}[section]

\newproclaim{rema}{Remark}[section]
\makeatother

\begin{document}
\begin{frontmatter}

\title{Limit theorems for power variations of pure-jump processes with application to activity estimation}

\runtitle{Limit theorems for power variations}

\begin{aug}
\author[A]{\fnms{Viktor} \snm{Todorov}\ead[label=e1]{v-todorov@northwestern.edu}\corref{}}
\and
\author[B]{\fnms{George} \snm{Tauchen}\ead[label=e2]{george.tauchen@duke.edu}}

\runauthor{V. Todorov and G. Tauchen}

\affiliation{Northwestern University and Duke University}
\address[A]{Department of Finance\\
Northwestern University\\
Evanston,
Illinois 60208-2001\\
USA\\
\printead{e1}} %adresu isvedimo komanda gale!
\address[B]{Department of Economics\\
Duke University\\
Durham, North Carolina 27708-0097\\
USA\\
\printead{e2}}
\end{aug}

% HISTORY:
\received{\smonth{5} \syear{2009}}
\revised{\smonth{3} \syear{2010}}

% ABSTRACT
\begin{abstract}
This paper derives the asymptotic behavior of realized power variation
of pure-jump It\^{o} semimartingales as the sampling frequency within a
fixed interval increases to infinity. We prove convergence in
probability and an associated central limit theorem for the realized
power variation as a function of its power. We apply the limit theorems
to propose an efficient adaptive estimator for the activity of
discretely-sampled It\^{o} semimartingale over a fixed interval.
\end{abstract}



\end{frontmatter}

%s1 ###
\section{Introduction}\label{sec:intro}

Realized power variation of a discretely sampled process can be defined
as the sum of the absolute values of the increments of the process
raised to a given power. The leading case is when the power is $2$,
which corresponds to the realized variance that is widely used in
finance. It is well known that under very weak conditions (see, e.g.,
\cite{JS}) the realized variance converges to the quadratic variation
of the process as the sampling frequency increases.  Powers other than
$2$ have also been used as a way to measure variation of the process
over a given interval in time as well as for estimation in parametric
or semiparametric settings. Recently, Ait-Sahalia and Jacod
\cite{SJ06} have used the realized power variation as a way to test for
presence of jumps on a given path and \cite{JT07} have used it to test
for common arrival of jumps in a multivariate context.

The limiting behavior of the realized power variation has been studied
in the continuous semimartingale case in \cite{BNS03} and
\cite{BNGJPS05}. Some of these results are extended by \cite{BNSW05} to
situations when jumps are present but only when they have no asymptotic
effect on the behavior of the realized power variation.  A
comprehensive study of the limiting behavior of the realized power
variation when the observed process is a continuous semimartingale plus
possible jumps is contained in \cite{J06b}. This work includes also
cases when jumps affect the limit of the realized power variation.

A common feature of the above cited papers is that the observed process
always contains a continuous martingale. At the same time there are
different applications, for example, for modeling internet traffic
\cite{TT07} or volume of trades \cite{ACD08} and asset volatility
\cite{TT08}, where pure-jump semimartingales, that is,
semimartingales without a continuous martingale and nontrivial
quadratic variation, seem to be more appropriate. Parametric models of
pure-jump type for financial prices and/or volatility have been
proposed in \cite{BNS01,COGARCH,CGMY03}, among
others. The main goal of this paper is to derive the limit behavior of
the realized power variation of \textit{pure-jump semimartingales}.

Some work has already been done in this direction. When the power
exceeds the (generalized) Blumenthar--Getoor index of the jump process,
it follows from \cite{Lepingle} and \cite{J06b} that the (unscaled)
realized power variation converges almost surely to the sum of jumps
raised to the corresponding power, which in general is not
\textit{predictable} (\cite{JS}, Definition I.2.1) although the exact
rate of this convergence is not known.

The limiting behavior of the realized power variation when the
\textit{power is less than the Blumenthal--Getoor index} is not known
in general (apart from the fact that it explodes). Here we concentrate
precisely on this case. We make an assumption of locally stable
behavior of the \Lvb measure of the jump process. That is we assume
that the \Lvb measure behaves like that of a stable process around
zero, while its behavior for the ``big'' jumps is left unrestricted.
This assumption allows us to derive the asymptotic behavior of the
realized power variation in this case. Unlike the case when the power
exceeds the Blumenthal--Getoor index, here the realized power variation
needs to be scaled down by a factor determined by the
Blumenthal--Getoor index and its limit is an integral of a
\textit{predictable} process. The latter is a direct measure for the
stochastic volatility of the discretely-observed process, which is of
key interest for financial applications. Thus the realized power
variation for powers less than the Blumenthal--Getoor index contains
information for the value of this index as well as the underlying
stochastic volatility, and hence the importance of the limit results
for this range of powers that are derived here. Finally, in earlier
work \cite{Woerner03,Woerner,Woerner07}, some limit theorems for
realized power variations for pure-jump processes were studied, but the
results apply in somewhat limiting situations regarding time-dependence
and presence of a drift term (i.e., an absolutely continuous process),
both of which are very important characteristics of financial data.

A distinctive feature of this paper is that the convergence results for
the realized power variation are derived on the space of functions of
the power equipped with the uniform topology. In contrast, all previous
work has characterized the limiting behavior for a fixed power. The
uniform convergence is important when one needs to use an infinite
number of powers in estimation or the power of the realized power
variation needs first to be estimated itself from the data. Such a case
is illustrated in an application of the limit theorems derived in the
paper.

Our application is for the estimation of the activity level of a
discretely observed process. The latter is the smallest power for which
the realized power variation does not explode (formally the infimum).
In the case of a pure-jump process the activity level is just the
Blumenthal--Getoor index of the jumps and when a continuous martingale
is present it takes its highest value of $2$. Apart from the importance
of the Blumenthal--Getoor index in itself, the activity level provides
information on the type of the underlying process (e.g., whether it
contains a continuous martingale or not). The latter determines the
appropriate scaling factor of the realized power variation in
estimating integrated volatility measures.

We use the realized power variation computed over two different
frequencies to estimate the activity level. The choice of the power is
critical as it affects both efficiency and robustness. We develop an
adaptive estimation strategy using our limit results. In a first step
we construct an initial consistent estimator of the activity, and then,
based on the first step estimator, we choose the optimal power to
estimate the activity on the second step.

The paper is organized as follows. Section~\ref{sec:setup} presents the
theoretical setup. Section~\ref{sec:limit} derives convergence in
probability and associated central limit theorems for the appropriately
scaled realized power variation. Section~\ref{sec:est} applies the
limit results of Section~\ref{sec:limit} to propose an efficient
adaptive estimator of the activity of a discretely sampled process.
Section~\ref{sec:numer} contains a short Monte Carlo study of the
behavior of the estimator. Proofs are given in Section~\ref{sec:proof}.

%s2 ###
\section{Theoretical setup}\label{sec:setup}
The theoretical setup of the paper is as follows. We will assume that
we have discrete observations of some one-dimensional process, which we
will always denote with $X$. The process will be defined on some
filtered probability space $(\Omega,\mathcal{F},\mathbb{P})$ with
$\mathbf{F}$ denoting the filtration. We will restrict attention to the
class of \Ito semimartingales, that is, semimartingales with absolutely
continuous characteristics (see, e.g., \cite{JS}).

Throughout we will fix the time interval to be $[0,T]$, and we will
suppose that we observe the process $X$ at the equidistant times
$0,\Delta_n,\ldots,[T/\Delta_n]\Delta_n$, where $\Delta_n>0$. The
asymptotic results in this paper will be of fill-in type, that is, we
will be interested in the case when $\Delta_n\downarrow 0$ for a fixed
$T>0$.

The \textit{activity of the jumps} in $X$ is measured by the so-called
(generalized) Blumenthal--Getoor index. All of our limiting results for
the realized power variation will depend in an essential way on it. The
index is defined as
%e2.1 ###
\begin{equation}\label{eq:BG}
\inf\biggl\{r>0\dvtx \sum_{0\leq s\leq T}|\Delta X_s|^r<\infty\biggr\},
\end{equation}
where $\Delta X_s:=X_s-X_{s-}$. The index was originally defined in
\cite{BG61} only for pure-jump \Lvb processes. The definition in
\refeq{eq:BG} extends it to an arbitrary jump semimartingale and was
proposed in \cite{SJ07}. We recall the following well-known facts: (1)~the index takes its values in $[0,2]$; (2) it depends on the particular
realization of the process on the given interval; (3) the value of $1$
for the index separates finite from infinite variation jump processes.

Finally, we define the main object of our study, the realized power
variation. It is constructed from the discrete observations of the
process as
%e2.2 ###
\begin{equation}\label{eq:rpv}
V_t(p,X,\Delta_n)=\sum_{i=1}^{[t/\Delta_n]}|\Delta_i^n X|^p,\qquad p>0,\ t>0,
\end{equation}
where $\Delta_i^n X:= X_{i\Delta_n}-X_{(i-1)\Delta_n}$. Our main focus
will be the behavior of $V_t(p,X,\Delta_n)$ when $X$ is pure-jump
semimartingale and we will restrict further attention to the case when
the power is below the Blumenthal--Getoor index and the drift term has
no asymptotic effect.

%s3 ###
\section{Limit theorems for power variation}\label{sec:limit}

We start with deriving the asymptotic limit of the appropriately scaled
realized power variation and then proceed with a central limit theorem
associated with it. To ease exposition we first present the results in
the \Lvb case and then generalize to the case when $X$ is a
semimartingales with time-varying characteristics. For completeness we
state corresponding results in the case when $X$ is a continuous
martingale (plus jumps) as well.

%s3.1 ###
\subsection{Convergence in probability results}\label{subsec:cp}
The convergence in probability results have been already derived in
\cite{BNS03,BNGJPS05,J06b,Woerner,Woerner03,TT07} among others with various degrees of generality. We
briefly summarize them here as a starting point of our analysis. We
first introduce some notation that will be used throughout. We set
$\mu_p(\beta):=\mathbb{E}(|Z|^{p})$, where $Z$ is a random variable
with a standard stable distribution with index $\beta$ if $\beta<2$
[i.e., with characteristic function
$\mathbb{E}(\exp(iuZ))=\exp(-|u|^{\beta})$], and with standard normal
distribution if $\beta=2$ (i.e., normal with mean 0 and variance 1).
Further,
$\mu_{p,q}(\beta):=\mathbb{E}|Z^{(1)}|^{p_1}|Z^{(1)}+Z^{(2)}|^{p_2}$,
where $Z^{(1)}$ and $Z^{(2)}$ are two independent random variables
whose distribution is standard stable with index $\beta$ if $\beta<2$
and is standard normal if $\beta=2$. Finally, we denote
$\Pi_{A,\beta}:=2A\int_0^{\infty}(\frac{1-\cos(x)}{x^{\beta+1}})\,dx$ for
$\beta\in (0,2)$ and $A>0$.

Throughout, $\kappa(x)$ will denote a continuous truncation function,
that is, a continuous function with bounded support such that
$\kappa(x)\equiv x$ around the origin, and $\kappa'(x):=x-\kappa(x)$.

%s3.1.1 ###
\subsubsection{The \Lvb case}\label{subsec:cp_levy}

\begin{thms}\label{thms:cp_levy}
\textup{(a)} Suppose $X$ is given by
%e3.1 ###
\begin{equation}\label{eq:X_c}
dX_t=m_c\,dt+\sigma\,dW_t+\int_{\mathbb{R}}\kappa(x)\widetilde{\mu}(dt,dx)+\int_{\mathbb{R}}\kappa'(x)\mu(dt,dx),
\end{equation}
where $m_c$ and $\sigma\neq 0$ are constants, and $W_t$ is a standard
Brownian motion; $\mu$ is a homogenous Poisson measure with compensator
$F(dx)\,dt$. Denote with $\beta'$ the Blumenthal--Getoor index of the
jumps in $X$. Then, if $\beta'<2$ and for a fixed $T>0$, we have
%e3.2 ###
\begin{equation}\label{eq:cp_levy_a}
\Delta_n^{1-p/2}V_T(X,p,\Delta_n)\stackrel{\mathbb{P}}{\longrightarrow}T|\sigma|^{p}\mu_p(2),
\end{equation}
locally uniformly in $p\in(0,2)$.

\textup{(b)} Suppose $X$ is given by
%e3.3 ###
\begin{equation}\label{eq:X_d}
dX_t=m_d\,dt+\int_{\mathbb{R}}\kappa(x)\widetilde{\mu}(dt,dx)+\int_{\mathbb{R}}\kappa'(x)\mu(dt,dx),
\end{equation}
where $m_d$ is some constant; $\mu$ is a Poisson measure with
compensator $\nu(x)\,dx$ where
%e3.4 ###
\begin{equation}\label{eq:nu}
\nu(x)=\nu_1(x)+\nu_2(x),
\end{equation}
with
%e3.5 ###
\begin{equation}\label{eq:nu_1_2}
\nu_1(x)=\frac{A}{|x|^{\beta+1}}\quad\mbox{and}\quad
|\nu_2(x)|\leq\frac{B}{|x|^{\beta^{'}+1}}\qquad\mbox{when }|x|\leq x_0
\end{equation}
for some $A>0$, $B\geq 0$ and $x_0>0$; $\beta\in (0,2)$ and
$\beta'<\beta$. Assume that $m_d-\int_{\mathbb{R}}\kappa(x)\nu(x)\,dx=0$
if $\beta\leq 1$. Then for a fixed $T>0$, we have
%e3.6 ###
\begin{equation}\label{eq:cp_levy_b}
\Delta_n^{1-p/\beta}V_T(X,p,\Delta_n)\stackrel{\mathbb{P}}{\longrightarrow}T\Pi_{A,\beta}^{p/\beta}\mu_p(\beta),
\end{equation}
locally uniformly in $p\in(0,\beta)$.
\end{thms}

\begin{rema}
The crucial assumption in the pure-jump
case is the decomposition of the \Lvb measure in \refeq{eq:nu}. This
assumption implies that locally the process behaves like the stable,
that is, the very small jumps of the process are as if from a stable
process. This assumption allows to scale the realized power variation
using the Blumenthal--Getoor index $\beta$. We note that $\nu_2(x)$ is
not necessarily a \Lvb measure (since it can be negative) and thus
\refeq{eq:nu_1_2} does not allow to represent $X$ (in distribution) as
a sum of two independent jump processes, the first being the stable and
the second with Blumenthal--Getoor index of $\beta'$.
\end{rema}

\begin{rema}
If jumps are of finite variation, in
part (b) of the theorem we restrict $X$ to be equal to the sum of the
jumps on the interval. The reason for this is that if a drift term is
present (or equivalently a compensator for the small jumps), then it
``dominates'' the jumps and determines the behavior of the realized
power variation (see, e.g., \cite{J06b}).
\end{rema}

\begin{rema}
When $p>\beta$ in the pure-jump case the
limit of the realized power variation is just the some of the $p$th
absolute power of the jumps, and this result does not follow from a law
of large numbers but rather by proving that an approximation error for
this sum vanishes almost surely. Thus the behavior of the realized
power variation for $p<\beta$ and $p>\beta$ is fundamentally different.
The case $p=\beta$ is the dividing one. In this case the realized power
variation (unscaled) converges neither to a constant nor to the sum of
the absolute values of the jumps raised to the power $\beta$ (which is
infinite). It can be shown that after subtracting the ``big''
increments, that is, keeping only those for which $|\Delta_i^n X|\leq
K\Delta_n^{1/\beta}$, for an arbitrary constant $K>0$, the realized
power variation converges to a nonrandom constant.

We note that the behavior of the realized power variation for
$p\geq \beta$ in the pure-jump case is very different from the case
when $X$ \textit{does not} contain jumps. In the latter case for all
powers $(p\lesseqgtr2)$ the limit of the realized power variation is
determined by law of large numbers, and hence we always need to scale
the realized power variation in order to converge to a nondegenerate
limit (see, e.g., \cite{BNGJPS05}).
\end{rema}

%s3.1.2 ###
\subsubsection{Extension to general semimartingales}\label{subsec:cp_gen}

Now we extend Theorem~\ref{thms:cp_levy} to the case when $\sigma$ and
$\nu$ (and the drift terms $m_c$ and $m_d$) in \refeq{eq:X_c} and
\refeq{eq:X_d} are stochastic. Nothing fundamentally changes, apart
from the fact that the limits are now random (depending on the
particular realization of the process $X$). In the case of continuous
martingale plus jumps, we can substitute \refeq{eq:X_c} with the
following:
\begin{eqnarray}\label{eq:X_c_sv}
dX_t
&=&
m_{ct}\,dt+\sigma_{1t}\,dW_t+\int_{\mathbb{R}}\kappa(\delta(t,x))\widetilde{\mu}(dt,dx)\nonumber
\\[-8pt]\\[-8pt]
&&{}+
\int_{\mathbb{R}}\kappa'(\delta(s,x))\mu(dt,dx),\nonumber
\end{eqnarray}
where $m_{ct}$ is locally bounded and $\sigma_{1t}$ is a process with
\cadlag paths; in addition $|\sigma_{1t}|>0$ and $|\sigma_{1t-}|>0$ for
every $t>0$ almost surely; $\mu$ is a homogenous Poisson measure with
compensator $F(dx)\,dt$ and $\delta(t,x)$ is a predictable function
satisfying
%e3.7 ###
\begin{eqnarray}\label{eq:X_c_sv_a}
&&\mbox{the process $t\rightarrow\sup_{x}\frac{|\delta(t,x)|}{\gamma(x)}$ is locally bounded with}\nonumber
\\
&&\mbox{$\int_{\mathbb{R}}(|\gamma(x)|^{\beta'}\wedge 1)F(dx)<\infty$ for some nonrandom function $\gamma(x)$}
\\
&&\mbox{and some constant $\beta'\in[0,2]$.}\nonumber
\end{eqnarray}
Additionally we assume that $\sigma_{1t}$ is an \Ito semimartingale
satisfying equations similar to \refeq{eq:X_c_sv} and \refeq{eq:X_c_sv_a}
(with arbitrary driving Brownian motion and Poisson measure (and jump
size function) satisfying a condition as \refeq{eq:X_c_sv_a} with
$\beta'=2$) with locally bounded coefficients. We note that the
generalized Blumenthal--Getoor index of the jumps of $X$ in
\refeq{eq:X_c_sv} is bounded by the nonrandom $\beta'$.

In the pure-jump case more care is needed in introducing time
variation. Essentially we should keep the behavior around $0$ of the
jump compensator intact. Therefore the generalization of \refeq{eq:X_d}
that we consider is given by
%e3.8 ###
\begin{equation}\label{eq:X_d_sv}
dX_t=m_{dt}\,dt+\int_{\mathbb{R}}\sigma_{2t-}\kappa(x)\widetilde{\mu}(dt,dx)+\int_{\mathbb{R}}\sigma_{2t-}\kappa'(x)\mu(dt,dx),
\end{equation}
where $m_{dt}$ and $\sigma_{2t}$ are processes with \cadlag paths;
$\mu$ is a jump measure with compensator $\nu(x)\,dx\,dt$ where $\nu(x)$ is
given by \refeq{eq:nu}. We note that under this specification, the
generalized Blumenthal--Getoor index of $X$ in \refeq{eq:X_d_sv} equals
$\beta$ on \textit{every path}, where $\beta$ is the constant appearing
in \refeq{eq:nu_1_2}. Further we assume $|\sigma_{2t}|>0$ and
$|\sigma_{2t-}|>0$ for every $t>0$ almost surely and impose the
following dynamics for the process $\sigma_{2t}$:
\begin{eqnarray}\label{eq:X_d_sv_sv}
d\sigma_{2t}
&=&
b_{2t}\,dt+\widetilde{\sigma}_{2t}\,dW_t+\int_{\mathbb{R}^2}\kappa(\delta(t,x))\widetilde{\underline{\mu}}(dt,d\mathbf{x})\nonumber
\\[-8pt]\\[-8pt]
&&{}+
\int_{\mathbb{R}^2}\kappa'(\delta(t,x))\underline{\mu}(dt,d\mathbf{x}),\nonumber
\end{eqnarray}
where $W$ is a Brownian motion; $\underline{\mu}$ is a homogenous
Poisson measure on $\mathbb{R}^2$ with compensator
$\underline{\nu}(d\mathbf{x})\,dt$ for $\underline{\nu}$ denoting some
$\sigma$-finite measure on $\mathbb{R}^2$, satisfying
$\underline{\mu}(dt,A\times \mathbb{R})\equiv \mu(dt,A)$ for any $A\in
\mathcal{B}(\mathbb{R}_0)$ with
$\mathbb{R}_0:=\mathbb{R}\setminus\{0\}$; $\delta(t,\mathbf{x})$ is an
$\mathbb{R}$-valued predictable function satisfying
%e3.9 ###
\begin{eqnarray}\label{eq:X_d_sv_a}
&&
\mbox{the process $t\rightarrow
\sup_{\mathbf{x}}\frac{|\delta(t,\mathbf{x})|}{\gamma(\mathbf{x})}$ is locally bounded with}\nonumber
\\
&&
\mbox{$\int_{\mathbb{R}^2}(|\gamma(\mathbf{x})|^{\beta+\varepsilon}\wedge 1)\underline{\nu}(\mathbf{x})\,d\mathbf{x}<\infty$ for some nonrandom function}
\\
&&
\mbox{on $\mathbb{R}^2$, $\gamma(\mathbf{x})$, where $\beta$ is the constant in \refeq{eq:nu_1_2}, and for $\forall \varepsilon>0$.}\nonumber
\end{eqnarray}
Additionally we assume that $m_{dt}$ and $\widetilde{\sigma}_{2t}$ are
\Ito semimartingales satisfying equations similar to
\refeq{eq:X_c_sv} and \refeq{eq:X_c_sv_a} (with arbitrary driving Brownian
motion and Poisson measure) with locally bounded coefficients. This
specification for $\sigma_{2t}$ is fairly general and it importantly
allows for dependence between the driving jump measure in
\refeq{eq:X_d_sv} and $\sigma_{2t}$, which is important for financial
applications (see, e.g., the COGARCH model of \cite{COGARCH}).

The restrictions on $\sigma_{1t}$ and $\sigma_{2t}$ in
\refeq{eq:X_c_sv} and \refeq{eq:X_d_sv_sv} are stronger than needed for
the convergence in probability results in the next theorem, but they
will be used for deriving the central limit results in the next
subsection. These assumptions are nevertheless weak and therefore we
impose them throughout. For example, the \Ito semimartingale
restrictions on $\sigma_{1t}$ and $\sigma_{2t}$ and their coefficients,
together with conditions \refeq{eq:X_c_sv_a} and \refeq{eq:X_d_sv_a},
will be automatically satisfied if $X$ solves
%e3.10 ###
\begin{equation}\label{eq:levy_sde}
dX_t=f(X_{t-})\,dL_t
\end{equation}
for some twice continuously differentiable function $f(\cdot)$ with at
most linear growth and $L$ being the \Lvb process in \refeq{eq:X_c} or
\refeq{eq:X_d} (see, e.g., Remark 2.1 in \cite{J06b}). The next theorem
states the general result on convergence in probability of realized
power variation.

\begin{thms}\label{thms:cp_gen}
\textup{(a)} Suppose $X$ is given by \refeq{eq:X_c_sv} and \refeq{eq:X_c_sv_a}
is satisfied with $\beta'<2$. Then for a fixed $T>0$ we have
%e3.11 ###
\begin{equation}\label{eq:cp_gen_a}
\Delta_n^{1-p/2}V_T(X,p,\Delta_n)\stackrel{\mathbb{P}}{\longrightarrow}\mu_p(2)\int_0^T|\sigma_{1s}|^{p}\,ds,
\end{equation}
locally uniformly in $p\in(0,2)$.

\textup{(b)} Suppose $X$ is given by \refeq{eq:X_d_sv}, \refeq{eq:X_d_sv_sv} and
\refeq{eq:nu_1_2} holds with $\beta'<\beta$. Further assume
$m_{ds}-\sigma_{2s-}\int_{\mathbb{R}}\kappa(x)\nu(x)\,dx$ is identically
zero on $[0,T]$ on the observed path if $\beta\leq 1$. Then for a fixed
$T>0$ we have
%e3.12 ###
\begin{equation}\label{eq:cp_gen_b}
\Delta_n^{1-p/\beta}V_T(X,p,\Delta_n)\stackrel{\mathbb{P}}{\longrightarrow}\Pi_{A,\beta}^{p/\beta}\mu_p(\beta)\int_0^T|\sigma_{2s}|^p\,ds,
\end{equation}
locally uniformly in $p\in(0,\beta)$.
\end{thms}

\begin{rema}
As seen from the above theorem, in both
cases the (scaled) realized power variation estimates an integrated
volatility measure $\int_0^T|\sigma_{is}|^p\,ds$ for $i=1,2$, which is
important for measuring volatility in financial applications. What is
different in the two cases is the scaling factor that is used. The
latter depends on the \textit{activity} of $X$ that we formally define
later in Section~\ref{sec:est} and then estimate using the limit
theorems of the current section.
\end{rema}

%s3.2 ###
\subsection{CLT results}\label{subsec:clt}

Since in our application we make use of the realized power variation
over two frequencies, $\Delta_n$ and $2\Delta_n$, we derive a CLT for
the vector $(V_{T}(X,p,2\Delta_n),V_{T}(X,p,\Delta_n))'$. In the next
and subsequent theorems $\mathcal{L}-s$ will stand for convergence
stable in law (see, e.g., \cite{JS} for a definition for filtered
probability spaces).

%s3.2.1 ###
\subsubsection{The \Lvb case}\label{subsec:clt_levy}

As for the convergence in probability we start with the \Lvb case. The
result is given in the following theorem.

\begin{thms}\label{thms:asl}
\textup{(a)} Suppose $X$ is given by the process in \refeq{eq:X_c} with
Blumen\-thal--Getoor index $\beta'<1$. Then, for a fixed $T>0$ and any
$0<p_l\leq p_h<1$ such that $\frac{\beta'}{2-\beta'}<p_l\leq p_h<1$, we
have
%e3.13 ###
\begin{equation}\label{eq:asl_a}
\qquad\Delta_n^{-1/2}
\pmatrix{
\Delta_n^{1-p/2}V_{T}(X,p,2\Delta_n)-2^{p/2-1}T|\sigma|^p\mu_p(2)\cr
\Delta_n^{1-p/2}V_{T}(X,p,\Delta_n)-T|\sigma|^p\mu_p(2)}
\stackrel{\mathcal{L}-s}{\longrightarrow}\Psi_{2,T}(p),
\end{equation}
where the convergence takes place in
$\mathcal{C}(\mathbb{R}^2,[p_l,p_h])$, the space of
$\mathbb{R}^2$-valued continuous functions on $[p_l,p_h]$ equipped with
the uniform topology; $\Psi_{2,T}(p)$ is a continuous centered Gaussian
process, independent from the filtration on which $X$ is defined, with
the following variance--covariance $\cov (\Psi_{2,T}(p),\allowbreak\Psi_{2,T}(q))$
for some $p,q\in [p_l,p_h]$:
\[
T|\sigma|^{2p}
\pmatrix{
2^{(p+q)/2-1}\bigl(\mu_{p+q}(2)-\mu_p(2)\mu_q(2)\bigr)&\mu_{q,p}(2)-2^{p/2}\mu_p(2)\mu_q(2)\cr
\mu_{p,q}(2)-2^{q/2}\mu_p(2)\mu_q(2) & \mu_{p+q}(2)-\mu_p(2)\mu_q(2)}.
\]

\textup{(b)} Suppose $X$ is given by the process in \refeq{eq:X_d}, and
\refeq{eq:nu_1_2} holds with $\beta'<\beta/2$. Then, for a fixed $T>0$
and any $0<p_l\leq p_h<1$ such that either \textup{(i)}
$(\frac{2-\beta}{2(\beta-1)}\vee
\frac{\beta\beta'}{2(\beta-\beta')})<p_l\leq p_h<\beta/2$ when
$\beta>\sqrt{2}$ or \textup{(ii)} $m_d\equiv 0$, $\nu$ and $\kappa$ symmetric
and $\frac{\beta\beta'}{2(\beta-\beta')}<p_l\leq p_h<\beta/2$, we have
\begin{eqnarray}\label{eq:asl_b}
&&
\Delta_n^{-1/2}
\pmatrix{
\Delta_n^{1-p/\beta}V_{T}(X,p,2\Delta_n)-2^{p/\beta-1}T\Pi_{A,\beta}^{p/\beta}\mu_p(\beta)\cr
\Delta_n^{1-p/\beta}V_{T}(X,p,\Delta_n)-T\Pi_{A,\beta}^{p/\beta}\mu_p(\beta)}\nonumber
\\[-8pt]\\[-8pt]
&&\qquad\stackrel{\mathcal{L}-s}{\longrightarrow}
\Psi_{\beta,T}(p),\nonumber
\end{eqnarray}
where the convergence takes place in
$\mathcal{C}(\mathbb{R}^2,[p_l,p_h])$, the space of
$\mathbb{R}^2$-valued continuous functions on $[p_l,p_h]$ equipped with
the uniform topology; $\Psi_{\beta,T}(p)$ is a continuous centered
Gaussian process, independent from the filtration on which $X$ is
defined, with the following variance--covariance $\cov
(\Psi_{\beta,T}(p),\allowbreak\Psi_{\beta,T}(q))$ for some $p,q\in [p_l,p_h]$:
\[
T\Pi_{A,\beta}^{2p/\beta}\!
\pmatrix{2^{(p+q)/\beta-1}\bigl(\mu_{p+q}(\beta)-\mu_p(\beta)\mu_q(\beta)\bigr)\!&\mu_{q,p}(\beta)-2^{p/\beta}\mu_p(\beta)\mu_q(\beta)\cr
\mu_{p,q}(\beta)-2^{q/\beta}\mu_p(\beta)\mu_q(\beta)\!&\mu_{p+q}(\beta)-\mu_p(\beta)\mu_q(\beta)}.
\]
\end{thms}

\begin{rema}
The result in part (a) for a
\textit{fixed} $p$ has been already shown (see, e.g., \cite{BNGJPS05}
and references therein). In the pure-jump case \refeq{eq:X_d}, the
result in \refeq{eq:asl_b} for a \textit{fixed} $p$ has been derived by
\cite{Woerner03} but only in the case \textit{when there is no drift}
[i.e., only under condition (ii) in part (b) of Theorem~\ref{thms:asl}]
and a slightly more restrictive condition on the residual measure
$\nu_2$. The general treatment here is important for financial
applications, as the presence of risk premium means theoretically that
the dynamics of traded assets should contain a drift term. Allowing for
a drift term is also important for applications to processes exhibiting
strong mean reversion like asset volatilities and trading volumes (see,
e.g., \cite{ACD08}).
\end{rema}

\begin{rema}
Theorem~\ref{thms:asl} shows that the
convergence of the scaled and centered power variation is
\textit{uniform over $p$}. This result has not been shown before. The
uniformity is important, for example, in adaptive estimation where the
power of the realized power variation to be used needs to be estimated
from the data. This is illustrated in our application in
Section~\ref{sec:est}.
\end{rema}

\begin{rema}
Comparing Theorem~\ref{thms:asl} with
Theorem~\ref{thms:asp} we see that both in parts (a) and (b) we have
imposed the stricter restrictions,
\[
p\in\biggl(\frac{2-\beta}{2(\beta-1)}\vee\frac{\beta\beta'}{2(\beta-\beta')},\beta/2\biggr)
\]
[with $\beta=2$ for part (a)] and $\beta'<\beta/2$. The lower bound for $p$ is determined
from the presence of a ``less active'' component in $X$. The
restriction $p>\frac{2-\beta}{2(\beta-1)}$ comes from the presence of a
drift term. We note that it is more restrictive the lower the $\beta$
is. In fact when $\beta\leq \sqrt{2}$, the presence of a drift term
will slow down the rate of convergence of the scaled power variation,
and therefore the limiting result in \refeq{eq:asl_b} will not hold. In
contrast for high values of $\beta$, $p>\frac{2-\beta}{2(\beta-1)}$ is
very weak and in the limiting case when $\beta=2$ [part (a) of the
theorem] it is never binding. We can interpret the restrictions
$p>\frac{\beta\beta'}{2(\beta-\beta')}$ and $\beta'<\beta/2$ similarly.
They come from the presence in $X$ of a less active jump component with
Blumenthal--Getoor index $\beta'$.

Also, the restriction $p<\beta/2$, which in particular implies that the
function $|x|^p$ is subadditive, is crucial for bounding the effect of
the ``residual'' jump components in $X$.
\end{rema}

\begin{rema}\label{rema38}
We can also derive a central limit
theorem when $p\in (\beta/2,\beta)$ (and when there are no ``residual''
jump components). In this case pure-continuous and pure-jump
martingales differ. While in the former case the rate of convergence
continuous to be $\sqrt{\Delta_n}$, in the latter the rate slows down.
The precise result is the following:

\textit{Suppose $X$ is symmetric stable plus a drift, that
is, the process in \textup{\refeq{eq:X_d}} with $\nu_2(x)\equiv 0$ and further
$m_d-\int_{\mathbb{R}}\kappa(x)\nu_1(x)\,dx\equiv 0$ when $\beta\leq 1$.
Set $a=m_d+\int_{\mathbb{R}}(x-\kappa(x))\nu_1(x)\,dx$ when $\beta>1$ and
$a=0$ when $\beta\leq 1$. Then for  a fixed $p\in (\beta/2\vee
\frac{1}{\beta}1_{\{\beta>1\cap a\neq0\}},\beta)$ we have}
%e3.14 ###
\begin{equation}\label{eq:sclt}
\Delta_n^{p/\beta-1}\bigl(\Delta_n^{1-p/\beta}V_{T}(X,p,\Delta_n)-T\Pi_{A,\beta}^{p/\beta}\mu_p(\beta)\bigr)\stackrel{\mathcal{L}}{\longrightarrow}
S_{T},
\end{equation}
\textit{where $S_t$ is pure-jump \Lvb process with \Lvb density
$1_{\{x>0\}}2\frac{A}{p}\frac{1}{x^{1+\beta/p}}$ and zero drift with
respect to the ``truncation'' function $\kappa(x)=x$. This is an
asymmetric stable process with index $\beta/p\in (1,2)$. }

As seen from \refeq{eq:sclt}, as we increase $p$ the rate of
convergence of the realized power variation slows down from
$\sqrt{\Delta_n}$ to $1$. Therefore this range of powers is less
attractive for estimation purposes. This will be further discussed in
Section~\ref{sec:est}.
\end{rema}

%s3.2.2 ###
\subsubsection{Extension to general semimartingales}\label{subsec:clt_gen}

We proceed with the analogue of Theorem~\ref{thms:asl} in the more
general setup of Section~\ref{subsec:cp_gen}. We state the case when
$\beta>\sqrt{2}$ only, since as seen from Theorem~\ref{thms:asl} and
Remark \ref{rema38}, the case $\beta\leq\sqrt{2}$ needs an assumption of zero
drift, and this limits its usefulness for financial applications where
the drift arises from the presence of risk premium.

\begin{thms}\label{thms:sv}
\textup{(a)} Suppose $X$ is given by \refeq{eq:X_c_sv}, and \refeq{eq:X_c_sv_a}
is satisfied with $\beta'<1$. Then, for a fixed $T>0$ and any
$0<p_l\leq p_h<1$ such that $\frac{\beta'}{2-\beta'}<p_l\leq p_h<1$, we
have
\begin{eqnarray}\label{eq:sv_a}
&&\Delta_n^{-1/2}
\pmatrix{
\displaystyle\Delta_n^{1-p/2}V_{T}(X,p,2\Delta_n)-2^{p/2-1}\mu_p(2)\int_0^T|\sigma_{1s}|^p\,ds\cr
\displaystyle\Delta_n^{1-p/2}V_{T}(X,p,\Delta_n)-\mu_p(2)\int_0^T|\sigma_{1s}|^p\,ds}\nonumber
\\[-8pt]\\[-8pt]
&&\qquad\stackrel{\mathcal{L}-s}{\longrightarrow}
\Psi_{2,T}(p),\nonumber
\end{eqnarray}
where the convergence takes place in
$\mathcal{C}(\mathbb{R}^2,[p_l,p_h])$, the space of
$\mathbb{R}^2$-valued continuous functions on $[p_l,p_h]$ equipped with
the uniform topology; $\Psi_{2,T}(p)$ is a continuous centered Gaussian
process, independent from the filtration on which $X$ is defined, with
the following variance--covariance $\cov (\Psi_{2,T}(p),\allowbreak\Psi_{2,T}(q))$
for some $p,q\in [p_l,p_h]$:
\begin{eqnarray*}
&&
\int_0^T|\sigma_{1s}|^{2p}\,ds
\\
&&\hphantom{\int_0^T}{}\times
\pmatrix{
2^{(p+q)/2-1}\bigl(\mu_{p+q}(2)-\mu_p(2)\mu_q(2)\bigr)&\mu_{q,p}(2)-2^{p/2}\mu_p(2)\mu_q(2)\cr
\mu_{p,q}(2)-2^{q/2}\mu_p(2)\mu_q(2)&\mu_{p+q}(2)-\mu_p(2)\mu_q(2)}.
\end{eqnarray*}

\textup{(b)} Suppose $X$ is given by \refeq{eq:X_d_sv}--\refeq{eq:X_d_sv_a} with
$\beta>\sqrt{2}$ and \refeq{eq:nu_1_2} holds with $\beta'<\beta/2$.
Then, for a fixed\vspace*{-2pt} $T>0$ and any $0<p_l\leq p_h<1$ such that
$(\frac{2-\beta}{2(\beta-1)}\vee \frac{\beta-1}{2}\vee\frac{\beta\beta'}{2(\beta-\beta')})<p_l\leq p_h<\beta/2$, we have
\begin{eqnarray}\label{eq:sv_b}
&&
\Delta_n^{-1/2}
\pmatrix{
\displaystyle\Delta_n^{1-p/\beta}V_{T}(X,p,2\Delta_n)-2^{p/\beta-1}\Pi_{A,\beta}^{p/\beta}\mu_p(\beta)\int_0^T|\sigma_{2s}|^p\,ds\cr
\displaystyle\Delta_n^{1-p/\beta}V_{T}(X,p,\Delta_n)-\Pi_{A,\beta}^{p/\beta}\mu_p(\beta)\int_0^T|\sigma_{2s}|^p\,ds}\nonumber
\\[-8pt]\\[-8pt]
&&\qquad\stackrel{\mathcal{L}-s}{\longrightarrow}
\Psi_{\beta,T}(p),\nonumber
\end{eqnarray}
where the convergence takes place in
$\mathcal{C}(\mathbb{R}^2,[p_l,p_h])$---the space of
$\mathbb{R}^2$-valued continuous functions on $[p_l,p_h]$ equipped with
the uniform topology; $\Psi_{\beta,T}(p)$ is a continuous centered
Gaussian process, independent from the filtration on which $X$ is
defined, with the following variance--covariance $\cov
(\Psi_{\beta,T}(p),\allowbreak\Psi_{\beta,T}(q))$ for some $p,q\in [p_l,p_h]$:
\begin{eqnarray*}
&&
\hspace{-5pt}\Pi_{A,\beta}^{2p/\beta}\int_0^T|\sigma_{2s}|^{2p}\,ds\nonumber
\\[-8pt]\\[-8pt]
&&\qquad{}\times
\pmatrix{
2^{(p+q)/\beta-1}\bigl(\mu_{p+q}(\beta)-\mu_p(\beta)\mu_q(\beta)\bigr)\!&\!\mu_{q,p}(\beta)-2^{p/\beta}\mu_p(\beta)\mu_q(\beta)\cr
\mu_{p,q}(\beta)-2^{q/\beta}\mu_p(\beta)\mu_q(\beta)\!&\!\mu_{p+q}(\beta)-\mu_p(\beta)\mu_q(\beta)}.\nonumber
\end{eqnarray*}
\end{thms}

Part (a) of the theorem has been derived in \cite{BNGJPS05}, while part
(b) is a new result. We note that compared with the \Lvb case in part
(b) of the theorem we have a slightly stronger restriction for $p$,
that is, $p$ cannot be arbitrarily small when $\beta$ is close to~2.
This is of no practical concern as the very low powers are not very
attractive because of the high associated asymptotic variance. This is
further discussed in Section~\ref{sec:est}.

%s4 ###
\section{Application: Adaptive estimation of activity}\label{sec:est}

We proceed with an application of our limit results. We first define
our object of interest, the activity level of the discretely-observed
process, and show how the realized power variation can be used for its
inference. Following that we develop an adaptive strategy for its
estimation.

%s4.1 ###
\subsection{Definitions}\label{subsec:est_def}

We define the \textit{activity level of an \Ito semimartingale} $X$ as
the smallest power for which the realized power variation does not
explode, that is,
%e4.1 ###
\begin{equation}\label{eq:activity}
\beta_{X,T}:= \inf\Bigl\{r>0\dvtx\operatorname{plim}\limits_{\displaystyle\Delta_n\rightarrow 0}V(r,X,\Delta_n)_T<\infty\Bigr\}.
\end{equation}
$\beta_{X,T}$ takes values in $[0,2]$ and is defined pathwise. It is
determined by the most active component in $X$ and the order of the
different components forming the \Ito semimartingale from least to most
active is: finite activity jumps, jumps of finite variation, drift
(absolutely continuous process), infinite variation jumps, continuous
martingale. When the dominating component of $X$ is its jump part (and
only then), $\beta_{X,T}$ coincides with the generalized
Blumenthal--Getoor index. Thus, for $X$ in \refeq{eq:X_c_sv},
$\beta_{X,T}\equiv 2$, and for $X$ in
\refeq{eq:X_d_sv} and \refeq{eq:X_d_sv_sv}, $\beta_{X,T}\equiv \beta$. We
note that $\beta_{X,T}$ determines uniquely the appropriate scale for
the realized power variation in the estimation of the integrated
volatility measures of the process (see
Theorems~\ref{thms:cp_levy} and \ref{thms:cp_gen}).

When the process is observed discretely, $\beta_{X,T}$ is unknown and
our goal is to derive an estimator for it. Since the scaling of the
realized power variation depends on the activity level, we can identify
the latter by taking a ratio of the realized power variation over two
scales. Therefore our estimation will be based on the following
function of the power:
%e4.2 ###
\begin{equation}\label{eq:asf}
\qquad b_{X,T}(p)=\frac{\ln(2)p}{\ln(2)+\ln[V_T(X,p,2\Delta_n)]-\ln[V_T(X,p,\Delta_n)]},\qquad p>0.
\end{equation}
A two-scale approach for related problems has been previously used also
in \cite{ZMA05,SJ07,TT07}.

%s4.2 ###
\subsection{Limit behavior of $b_{X,T}(p)$}\label{subsec:est_asf}

For ease of exposition here we restrict attention to the \Lvb case. The
extension to the general semimartingales in \refeq{eq:X_c_sv}, \refeq{eq:X_d_sv} and \refeq{eq:X_d_sv_sv} follows from an easy application
of Theorem~\ref{thms:sv}. In what follows, for any $p$ and $q$ both in
$(0,\beta/2)$ we denote
\begin{eqnarray}\label{eq:Kp}
\qquad K_{p,q}(\beta)
&=&
\frac{\beta^4}{\ln^2(2)pq\mu_{p}(\beta)\mu_{q}(\beta)}\bigl(3\mu_{p+q}(\beta)+\mu_{p}(\beta)\mu_{q}(\beta)\nonumber
\\[-8pt]\\[-8pt]
&&\hphantom{\frac{\beta^4}{\ln^2(2)pq\mu_{p}(\beta)\mu_{q}(\beta)}(}{}-
2^{1-p/\beta}\mu_{p,q}(\beta)-2^{1-q/\beta}\mu_{q,p}(\beta)\bigr).\nonumber
\end{eqnarray}

\begin{cor}\label{thms:asp}
\textup{(a)} Suppose $X$ is given by \refeq{eq:X_c}. Then for a fixed $T>0$ and
any $0<p_l\leq p_h<1$ we have
%e4.3 ###
\begin{equation}\label{eq:asp_a}
\sqrt{\frac{T}{\Delta_n}}\bigl(b_{X,T}(p)-2\bigr)\stackrel{\mathcal{L}-s}{\longrightarrow}Z_{2}(p),\qquad
\mbox{uniformly on }[p_l,p_h],
\end{equation}
where $Z_{2}(p)$ is a centered Gaussian process on $[p_l,p_h]$ with
$\cov (Z_{2}(p),Z_{2}(q))$ $=$ $K_{p,q}(2)$ for some $p,q\in [p_l,p_h]$
and independent from the filtration on which $X$ is defined, provided
$\beta'<1$ and $\frac{\beta'}{2-\beta'}<p_l\leq p_h<1$, where $\beta'$
is the Blumenthal--Getoor index of $X$.

\textup{(b)} Suppose $X$ is given by \refeq{eq:X_d}. Then for a fixed $T>0$ and
any $0<p_l\leq p_h<1$ we have
%e4.4 ###
\begin{equation}\label{eq:asp_b}
\sqrt{\frac{T}{\Delta_n}}\bigl(b_{X,T}(p)-\beta\bigr)\stackrel{\mathcal{L}-s}{\longrightarrow}Z_{\beta}(p),\qquad
\mbox{uniformly on }[p_l,p_h],
\end{equation}
where $Z_{\beta}(p)$ is a centered Gaussian process on $[p_l,p_h]$ with
$\cov (Z_{\beta}(p),Z_{\beta}(q))$ $=$ $K_{p,q}(\beta)$ for some
$p,q\in [p_l,p_h]$ and independent from the filtration on which $X$ is
defined, provided \refeq{eq:nu_1_2} holds with $\beta^{'}<\beta/2$ and
either \textup{(i)} $(\frac{2-\beta}{2(\beta-1)}\vee
\frac{\beta\beta'}{2(\beta-\beta')})<p_l\leq p_h<\beta/2$ when
$\beta>\sqrt{2}$ or \textup{(ii)} $m_d\equiv 0$, $\nu$ symmetric and
$\frac{\beta\beta'}{2(\beta-\beta')}<p_l\leq p_h<\beta/2$.
\end{cor}

As seen from the corollary, $b_{X,T}(p)$ will estimate the activity
level only for powers that are below the activity level, which of
course is unknown. Corollary~\ref{thms:asp} shows further that the power
is also crucial for the rate at which the activity level is estimated.
The range of values of $p$ for which $b_{X,T}(p)$ is
$\sqrt{\Delta_n}$-consistent for $\beta_{X,T}$ defined in
\refeq{eq:activity} depends on the activity of the most active part of
the process, but also on the activity of the less active parts, that
is, $\beta'$ in part (a) and $\beta'\vee 1$ in part (b). For example,
when the observed process is a continuous martingale plus jumps [part
(a) of the corollary], then the activity of the jumps needs to be
sufficiently low in order to estimate $\beta_{X,T}$ at a rate
$\sqrt{\Delta_n}$. Similar observation holds for the pure-jump case as
well. The activity of the less active components of $X$ is unknown but
we want an estimator of $\beta_{X,T}$ that is robust, in the sense that
it has $\sqrt{\Delta_n}$ rate of convergence for most values of
$\beta'$. Based on the corollary, this means that we need to use values
of $p$ that are ``sufficiently'' close to half of the activity level
$\beta_{X,T}/2$.

%f1 ###
\begin{figure}

\includegraphics{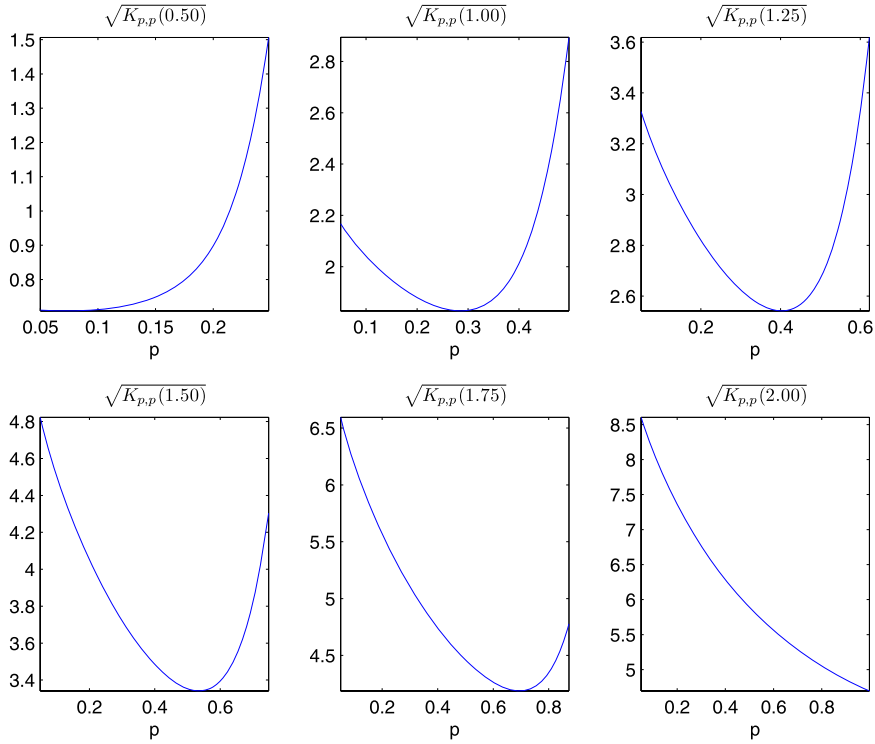}

\caption{Asymptotic
Standard Deviation of $b_{X,T}(p)$ for different values of $p$ and the
activity level $\beta_{X,T}$ defined in (\protect\ref{eq:activity}).
$K_{p,q}(\beta)$ is defined in (\protect\ref{eq:Kp}).}\label{fig:K}
\end{figure}

\begin{figure}

\includegraphics{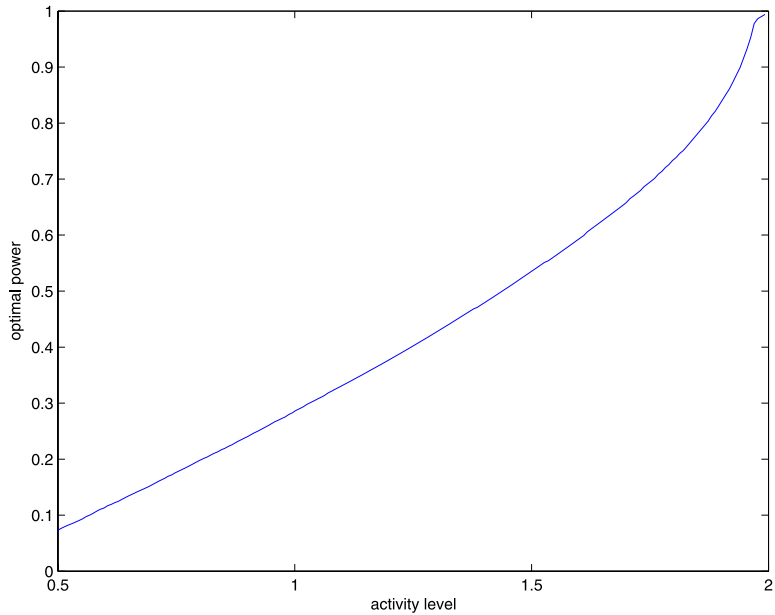}

\caption{Minimizing power $p$ of the asymptotic variance $K_{p,p}(\beta_{X,T})$ as a
function of the activity level $\beta_{X,T}$ defined in
(\protect\ref{eq:activity}).}\label{fig:asm}
\end{figure}

The presence of a less active component in the observed process aside,
the power at which $b_{X,T}(p)$ is evaluated is also important for the
rate of convergence and the asymptotic variance of the estimation of
the overall activity index. There is a difference between case (a) and
case (b) in this regard. When the activity level of $X$ is $2$ (and
there are no jumps), $b_{X,T}(p)$ will be $\sqrt{\Delta_n}$-consistent
for any power. In contrast, in the pure-jump case, this will be true
only for powers less than $\beta/2$.  Using powers
$p\in(\beta/2,\beta)$ slows down the rate of convergence from
$\sqrt{\Delta_n}$ to $1$, as pointed out in Remark \ref{rema38}. In
Figure~\ref{fig:K} we plotted the asymptotic standard deviation of
$b_{X,T}(p)$ for different values of the activity index $\beta_{X,T}$.
For activity less than $2$ the asymptotic variance has a pronounced
U-shape pattern, and as a result it is minimized somewhere within the
admissible range (for $\sqrt{\Delta_n}$-rate of convergence), but
\textit{the minimizing power depends on $\beta$}. On the other hand,
when $\beta_{X,T}=2$, i.e. when continuous martingale is present, the
asymptotic variance is minimized for  $p=1$ ($p=\beta_{X,T}/2$ is
admissible if $\beta_{X,T} \equiv 2$), although $\sqrt{K_{p,p}(2)}$
changes very little around $1$. These observations are further
confirmed from Figure~\ref{fig:asm}, which plots the power at which the
asymptotic variance is minimized as a function of the activity level.

%f2 ###

\begin{rema}
We note that in Corollary~\ref{thms:asp}
(and in fact throughout the paper) we kept $T$ fixed. What happens if
$T$ goes to infinity? In this case the result in
Corollary~\ref{thms:asp} will remain valid without any assumption on the
relative speed of $T\uparrow\infty$ and $\Delta_n\downarrow 0$ but only
in the case when $X$ is symmetric stable. In all other cases captured
by the specification in \refeq{eq:X_d} we will need to impose a
restriction on the relative speed with which $T$ increases. This
happens because the error in estimating $\beta_{X,T}$ depends on
$\Delta_n$ and \textit{cannot} vanish by just increasing the time span~$T$.
\end{rema}

%s4.3 ###
\subsection{Two-step estimation of activity}\label{subsec:est_tse}

We turn now to the explicit construction of an estimator of the
activity level guided by the results of Corollary~\ref{thms:asp}. Our
goal here is to derive a point estimator of the activity level which
has good robustness and efficiency properties. As we noted in the
previous subsection, the powers used in the construction of an
estimator for the activity level are crucial for its consistency, rate
of convergence and asymptotic efficiency. Importantly, whether to use a
given power in the estimation depends on the value of $\beta_{X,T}$
which is unknown and is itself being estimated.

This suggests implementing an adaptive (two-stage) estimation
procedure, where on a first stage we construct an initial
\textit{consistent} estimator of the activity. Any estimator with
arbitrary rate of convergence on this first stage can be used; the only
requirement is that it is consistent. Then, on a second stage, we can
use the first-stage estimator to select the power(s) at which
$b_{X,T}(p)$ is evaluated. This can be done because the convergence in
\refeq{eq:asp_a} and \refeq{eq:asp_b} is uniform in $p$. We give the
generic construction of the two-stage estimator in the \Lvb case in the
following theorem.

\begin{thms}\label{thms:ase}
Fix some $T>0$ and suppose $X$ is given either by \refeq{eq:X_c} or
\refeq{eq:X_d} with activity level $\beta_{X,T}$ defined in
\refeq{eq:activity}. Let $\hat{\beta}^{fs}_{X,T}$ be an arbitrary
consistent estimator of $\beta_{X,T}$ constructed from $X_0,
X_{\Delta_n},\ldots, X_{\Delta_n[T/\Delta_n]}$, that is, we have
$\hat{\beta}^{fs}_{X,T}\stackrel{\mathbb{P}}{\longrightarrow}\beta_{X,T}$
as $\Delta_n\rightarrow 0$. Suppose the functions $f_l(z)$ and
$f_h(z)$ are continuously differentiable in $z$ in a neighborhood of
$\beta_{X,T}$ and we have identically $0<f_l(z)<f_h(z)$. Set
\[
\tau_1^{*}=f_l(\beta_{X,T})\quad\mbox{and}\quad\tau_2^{*}=f_h(\beta_{X,T}),
\]
\[
\widehat{\tau}_1=f_l(\hat{\beta}^{fs}_{X,T})\quad\mbox{and}\quad\widehat{\tau}_2=f_h(\hat{\beta}^{fs}_{X,T}).
\]
Finally, denote
%e4.5 ###
\begin{equation}\label{eq:ts}
\hat{\beta}_{X,T}^{ts}=\int_{\widehat{\tau}_1}^{\widehat{\tau}_2}w(u)b_{X,t}(u)\,du,
\end{equation}
where $w(\cdot)$ is some weighting function, which is either continuous
on $[\tau_1^{*},\tau_2^{*}]$ or Dirac mass at some point in
$[\tau_1^{*},\tau_2^{*}]$ and such that
$\int_{\tau_1^{*}}^{\tau_2^{*}}w(u)\,du=1$. Then we have
%e4.6 ###
\begin{equation}\label{eq:ase_b}
\qquad\sqrt{\frac{T}{\Delta_n}}(\hat{\beta}_{X,T}^{ts}-\beta_{X,T})
\stackrel{\mathcal{L}-s}{\longrightarrow}
\varepsilon\times\sqrt{\int_{\tau_1^{*}}^{\tau_2^{*}}\!\!\int_{\tau_1^{*}}^{\tau_2^{*}}K_{u,v}(\beta_{X,T})w(u)w(v)\,du\,dv},
\end{equation}
where $\varepsilon$ is standard normal defined on an extension of the
original probability space provided:
\begin{enumerate}[(b)]
\item[(a)] if $X$ is given by \refeq{eq:X_c}, then
$\tau_2^{*}<\beta_{X,T}/2$ and the Blumenthal--Getoor index of the
jumps in $X$, $\beta'$, is such that
$\frac{\beta'}{2-\beta'}<\tau_1^{*}$ (which implies $\beta'<1$),
\item[(b)] if $X$ is given by \refeq{eq:X_d}, then $\tau_2^{*}<\beta/2$ and
either \textup{(i)} $\beta>\sqrt{2}$ and
$\tau_1^{*}>(\frac{2-\beta}{2(\beta-1)}\vee
\frac{\beta\beta'}{2(\beta-\beta')})$ or \textup{(ii)} $m_d\equiv 0$, $\nu$ and
$\kappa$ symmetric and
$\tau_1^{*}>\frac{\beta\beta'}{2(\beta-\beta')}$, where $\beta'$ is a
constant satisfying \refeq{eq:nu_1_2}.
\end{enumerate}
\end{thms}

The two-step estimator can be viewed as a weighted average of
$b_{X,T}(p)$ over an adaptively selected region of powers. This range
is determined on the basis of an initial consistent estimator of the
activity. The averaging of the powers on the second stage might be
beneficial since the correlation between the centered $b_{X,T}(p)$
evaluated over different powers is not perfect. We would expect that
the biggest benefit from averaging different powers in the estimation
will come from using powers that are sufficiently apart. However, as we
saw from Figure~\ref{fig:K}, significantly different powers would imply
that at least one of them is associated with too high asymptotic
variance and this could offset the benefit from the averaging.
Therefore, in practice on the second stage one can just evaluate
$b_{X,T}(p)$ at a single power. This case is stated in the next
corollary.

\begin{cor}\label{thms:tse}
Let $\hat{\beta}^{fs}_{X,T}$ be an arbitrary consistent\vspace*{-2pt} estimator
of $\beta_{X,T}$ constructed from $X_0, X_{\Delta_n},\ldots,X_{\Delta_n[T/\Delta_n]}$, that is, we have
$\hat{\beta}^{fs}_{X,T}\stackrel{\mathbb{P}}{\longrightarrow}\beta_{X,T}$
as $\Delta_n\rightarrow 0$. Set
%e4.7 ###
\begin{equation}\label{eq:est_ts}
\hat{\beta}^{ts}_{X,T}\equiv b_{X,T}(\widehat{\tau})\qquad\mbox{with }\widehat{\tau}:=f(\hat{\beta}^{fs}_{X,T}),
\end{equation}
where $f(\cdot)$ is some continuous function and further we set
$\tau^{*}:=f(\beta_{X,T})$. Then we have for a fixed $T$
%e4.8 ###
\begin{equation}\label{eq:ts_cor}
\sqrt{\frac{T}{\Delta_n}}(\hat{\beta}_{X,T}^{ts}-\beta_{X,T})\stackrel{\mathcal{L}-s}{\longrightarrow}\varepsilon\times\sqrt{K_{\tau^{*},\tau^{*}}(\beta_{X,T})}
\end{equation}
for $\varepsilon$ being standard normal, provided $\beta_{X,T}>2\tau^{*}$
and for $\beta'$ as in Theorem~\ref{thms:ase} we have:
\begin{enumerate}[(b)]
\item [(a)] if $X$ is given by \refeq{eq:X_c}, then $\beta'<\frac{2\tau^{*}}{1+\tau^{*}}$,
\item [(b)] if $X$ is given by \refeq{eq:X_d}, then $\beta'<\frac{2\beta\tau^{*}}{\beta+2\tau^{*}}$
and if $m_d\neq 0$ and/or $\nu$ is not symmetric then in addition we
also have $\beta>\sqrt{2}$ and $\tau^{*}<\frac{2-\beta}{2(\beta-1)}$.
\end{enumerate}
\end{cor}

A natural choice for the function $f(\cdot)$, that is, the power that
is used on the second stage, will be the one that minimizes the
asymptotic variance $K_{p,p}(\beta)$. This is further discussed in the
numerical implementation in the next section. Alternatively, one can
sacrifice some of the efficiency in exchange for robustness to a wider
range of $\beta'$ by picking power closer to $\beta_{X,T}/2$. We finish
this section with stating the equivalent of Corollary~\ref{thms:tse} in
the case when $X$ is a semimartingale with time-varying
characteristics. The theorem gives also feasible estimates of the
asymptotic variance of the two-step estimator.

\begin{thms}\label{thms:sv_est}
Suppose $\hat{\beta}_{X,T}^{fs}$ and $\hat{\beta}_{X,T}^{ts}$
are given by \refeq{eq:est_ts} for some fixed $T>0$.

\textup{(a)} If $X$ is given by \refeq{eq:X_c_sv} and \refeq{eq:X_c_sv_a} is
satisfied with $\beta'<\frac{2\tau^*}{1+\tau^*}$, then we have
%e4.9 ###
\begin{equation}\label{eq:ts_sv_a}
\frac{1}{\sqrt{\Delta_n}}(\hat{\beta}_{X,T}^{ts}-2)\stackrel{\mathcal{L}-s}{\longrightarrow}\varepsilon\times\sqrt{K_{\tau^*,\tau^*}(2)}\frac{\sqrt{\int_0^T|\sigma_{1s}|^{2\tau^*}\,ds}}{\int_0^T|\sigma_{1s}|^{\tau^*}\,ds},
\end{equation}
where $\varepsilon$ is standard normal and is defined on an extension of the original probability space.

\textup{(b)} If $X$ is given by \refeq{eq:X_d_sv}--\refeq{eq:X_d_sv_a} with
$\beta>\sqrt{2}$ and \refeq{eq:nu_1_2} holds with
$\beta'<\frac{\beta\tau^*}{1+\tau^*}$ and
$\tau^*\in(\frac{2-\beta}{2(\beta-1)}\vee \frac{\beta-1}{2},\beta/2)$,
then we have
%e4.10 ###
\begin{equation}\label{eq:ts_sv_b}
\frac{1}{\sqrt{\Delta_n}}(\hat{\beta}_{X,T}^{ts}-\beta)\stackrel{\mathcal{L}-s}{\longrightarrow}\varepsilon\times\sqrt{K_{\tau^*,\tau^*}(\beta)}\frac{\sqrt{\int_0^T|\sigma_{2s}|^{2\tau^*}\,ds}}{\int_0^T|\sigma_{2s}|^{\tau^*}\,ds},
\end{equation}
where $\varepsilon$ is standard normal and is defined on an extension of
the original probability space.

\textup{(c)} A consistent estimator for the asymptotic variance of both \refeq{eq:ts_sv_a} and \refeq{eq:ts_sv_b} is given by
%e4.11 ###
\begin{equation}\label{eq:ts_sv_c}
\Delta_n^{-1}K_{f(\hat{\beta}_{X,T}^{ts}),f(\hat{\beta}_{X,T}^{ts})}(\hat{\beta}_{X,T}^{ts})
\frac{\mu_{f(\hat{\beta}_{X,T}^{ts})}^2(\hat{\beta}_{X,T}^{ts})}{\mu_{2f(\hat{\beta}_{X,T}^{ts})}(\hat{\beta}_{X,T}^{ts})}
\frac{V_T(X,2f(\hat{\beta}_{X,T}^{ts}),\Delta_n)}{V_T^2(X,f(\hat{\beta}_{X,T}^{ts}),\Delta_n)}.
\end{equation}
\end{thms}

\begin{rema}
Although the choice of the first-step
estimator does not affect the first-order asymptotic properties of the
two-stage estimator, in practice it can matter a lot. One possible
choice for a first-step estimator of the activity is
%e4.12 ###
\begin{equation}
\widetilde{\beta}_{X,T}=\frac{2}{\ln(k)}\bigl(\ln(V'_T(\alpha,X,\Delta_n))-\ln(V'_T(\alpha,X,2\Delta_n))\bigr),
\end{equation}
where
$V'_T(\alpha,X,\Delta_n)=\sum_{i=1}^{[T/\Delta_n]}1_{\{|\Delta_i^nX|\geq
\alpha\sqrt{\Delta_n}\}}$ and $\alpha>0$ is an arbitrary constant. It
is easy to show that under the assumptions of Theorem~\ref{thms:sv},
$\widetilde{\beta}_{X,T}$ is a consistent estimator for $\beta_{X,T}$.
Another alternative first step estimator is $b_{X,T}(p)$ evaluated at
some small power. The latter will be a consistent estimator only if we
know apriori that the true value of $\beta_{X,T}$ is higher than some
positive number.
\end{rema}

%s5 ###
\section{Numerical implementation}\label{sec:numer}

In this section we test on simulated data the limit results of
Section~\ref{sec:limit}. We do this by investigating the finite sample
performance of the activity estimator of Section~\ref{sec:est}. In our
Monte Carlo study we work with the following model for $X$:
%e5.1 ###
\begin{equation}\label{eq:mc_x1}
X_t=\sigma_1W_t+\sigma_2\sum_{0\leq s\leq t}\Delta X_s,
\end{equation}
where the jumps of $X$ are with either of the following two
compensators:
%e5.2 ###
\begin{equation}\label{eq:mc_x2}
A\frac{e^{-\lambda |x|}}{|x|^{\beta+1}}\,dx\,ds\quad\mbox{or}\quad\lambda_{c}\delta_{\{x=\pm r\}}\,dx\,ds.
\end{equation}
The first compensator is that of a tempered stable \cite{CGMY02,Ro04} whose Blumenthal--Getoor index is the parameter $\beta$
and the second compensator is of a compound Poisson (which has of
course a Blumenthal--Getoor index of $0$). Note that for the tempered
stable process the value of $\beta'$ in \refeq{eq:nu_1_2} is equal to
$\beta-1\vee 0$. Therefore, the assumption $\beta'<\beta/2$ in
Theorems~\ref{thms:asl} and \ref{thms:sv} will always be satisfied.

%t1 ###
\begin{table}[b]
\caption{Parameter setting for the Monte Carlo}\label{tb:ps}
\begin{tabular}{@{}lccl@{}}
\hline
\textbf{Case}&$\bolds{\sigma_1^2}$&$\bolds{\sigma_2^2}$&\textbf{Jump specification}\\
\hline
A&0.0&1.0&Tempered stable with $A=1$, $\beta=1.50$ and $\lambda=0.25$\\
B&0.0&1.0&Tempered stable with $A=1$, $\beta=1.75$ and $\lambda=0.25$\\
C&0.8&0.0&None\\
D&0.8&1.0&Rare-jump with $\lambda_{c}=0.3333$, $r=0.7746$\\
\hline
\end{tabular}
\end{table}

In Table~\ref{tb:ps} we listed the four different cases we consider in
the Monte Carlo. The first two correspond to pure-jump processes with
two different values of the level of activity. The last two cases
correspond to a setting where a Brownian motion is present and
therefore overall activity of $X$ is $2$. In Case D the jumps in
addition to the Brownian motion have $20\%$ share in the total
variation of $X$ on a given interval, which is consistent with
empirical findings for financial price data.

If we think of a unit of time being a day, then in our Monte Carlo on
each ``day'' we sample $M=390$ times. This corresponds to approximately
every minute for $6.5$ hours trading day and every $5$ minutes for $24$
hours trading day. The activity estimation is performed over $22$ days,
that is, we set $T=22$. This corresponds to $1$ calendar month of
financial data. This Monte Carlo setup is representative of a typical
financial application that we have in mind. We do not report results
for other choices of $T$ and $M$ although we experimented with.  Quite
intuitively, an increase $T$ led to a reduction in the variance of the
estimators, while an increase in $M$ led to the elimination of any
existing biases. Finally, we consider $10{,}000$ number of Monte Carlo
replications.

Following our discussion in Section~\ref{subsec:est_tse} we calculate
over each simulation the following two-step estimator
$\hat{\beta}_{X,T}^{ts}$. In the first stage we evaluate the
function $b_{X,T}(p)$ at $p=0.1$. This yields an initial consistent,
albeit far from efficient, estimator for the activity, provided of
course the activity is above $0.1$. Then, given our first step
estimator of the activity, we compute the power at which
$K_{p,p}(\hat{\beta}^{fs}_{X,T})$ is minimized [recall the
definition of $K_{p,q}(\beta)$ in \refeq{eq:Kp}]. Our two-stage
estimator is simply the value of $b_{X,T}(p)$ at this optimal power.

%f3 ###
\begin{figure}[b]

\includegraphics{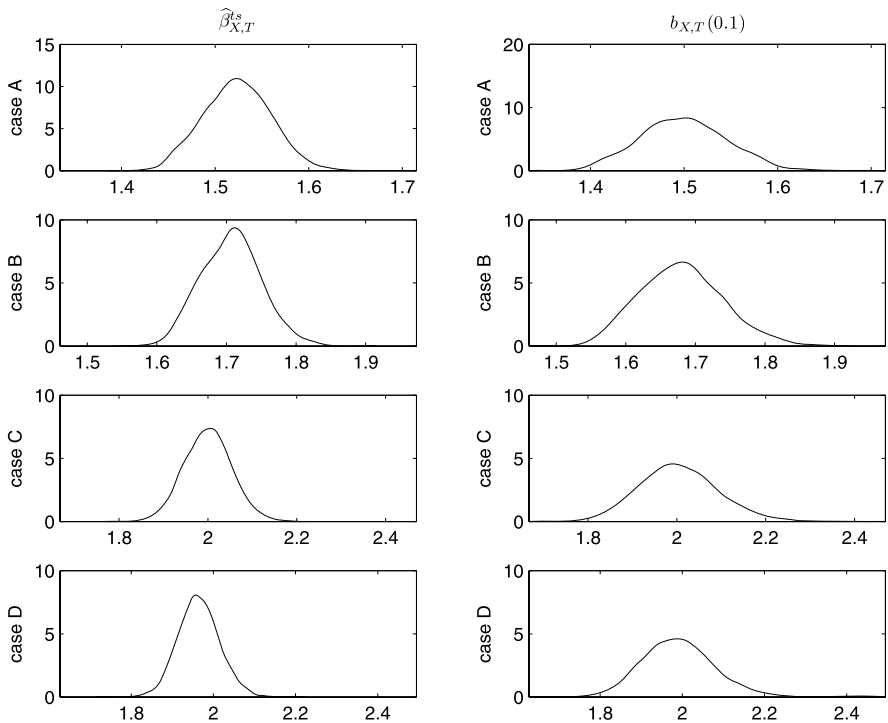}
 \caption{Histograms of $\hat{\beta}_{X,T}^{ts}$
and $b_{X,T}(0.1)$ from the Monte Carlo.}\label{fig:hist}
\end{figure}

In the Monte Carlo we compare the performance of our estimator with an
ad-hoc one where we simply evaluate $b_{X,T}(p)$ at the \textit{fixed}
``low'' power $p=0.1$. In Figure~\ref{fig:hist} we plot the histograms
of the two estimators $\hat{\beta}^{ts}_{X,T}$ and $b_{X,T}(0.1)$.
As we can see from this figure, the adaptive estimation of the activity
clearly outperforms the ad-hoc one based on a fixed power. In all cases
$\hat{\beta}_{X,T}^{ts}$ is much more concentrated around the true
value. This is further confirmed from Table~\ref{tb:comp}, which
reports summary statistics for the two estimators. The interquartile
range for the ad-hoc estimator is from $30\%$ to $60\%$ wider than that
of the adaptive estimator. A similar conclusion holds also for the mean
absolute deviation reported in the last column of the table. Thus, we
can conclude that choosing an ``optimal'' power can lead to nontrivial
improvements in the estimation of the activity, which is consistent
with our theoretical findings in Section~\ref{subsec:est_asf}.

%t2 ###
\begin{table}
\tablewidth=240pt
\caption{Comparison between two-step and one-step estimator}\label{tb:comp}
\begin{tabular*}{240pt}{@{\extracolsep{\fill}}lcccc@{}}
\hline
&\multicolumn{4}{c@{}}{\textbf{Summary statistics}}\\[-5pt]
&\multicolumn{4}{c@{}}{\hrulefill}\\
\textbf{Estimator}&$\bolds{\beta}$&\textbf{Median}&\textbf{IQR}&\textbf{MAD}\\
\hline
\multicolumn{5}{@{}l@{}}{\textit{Case A }}\\
\quad$\beta_{X,T}^{ts}$&1.50&1.5237&0.0495&0.0247\\
\quad$b_{X,T}(0.1)$&1.50&1.4985&0.0632&0.0316\\
\multicolumn{5}{@{}l@{}}{\textit{Case B }}\\
\quad$\beta_{X,T}^{ts}$&1.75&1.7075&0.0590&0.0294\\
\quad$b_{X,T}(0.1)$&1.75&1.6785&0.0814&0.0407\\
\multicolumn{5}{@{}l@{}}{\textit{Case C }}\\
\quad$\beta_{X,T}^{ts}$&2.00&2.0001&0.0719&0.0359\\
\quad$b_{X,T}(0.1)$&2.00&2.0005&0.1176&0.0588\\
\multicolumn{5}{@{}l@{}}{\textit{Case D }}\\
\quad$\beta_{X,T}^{ts}$&2.00&1.9632&0.0664&0.0332\\
\quad$b_{X,T}(0.1)$&2.00&1.9865&0.1164&0.0573\\
\hline
\end{tabular*}
\tabnotetext[]{tt1}{\textit{Note:} IQR is the inter-quartile range, and MAD is the mean
absolute deviation.}
\end{table}

We next investigate how well we can apply the feasible CLT for the
two-step activity estimator. For each estimated
$\hat{\beta}_{X,T}^{ts}$ we calculate standard errors using
\refeq{eq:ts_sv_c}. Table~\ref{tb:ase} provides summary statistics for
how well these estimated asymptotic standard errors track the exact
finite-sample standard error of the two-step estimator
$\hat{\beta}_{X,T}^{ts}$. Since $X$ is simulated from a \Lvb
process, the latter is computed as the standard error of
$\hat{\beta}_{X,T}^{ts}$ over the Monte Carlo replications.

%t3 ###
\begin{table}
\tablewidth=245pt
\caption{Precision of standard error estimation for the two-step
estimator}\label{tb:ase}
\begin{tabular*}{245pt}{@{\extracolsep{\fill}}lccc@{}}
\hline
&\multicolumn{3}{c@{}}{\textbf{Summary statistics for} $\bolds{\widehat{\mathbf{Ase}}}\bolds{(\beta_{X,T}^{ts}))}$}\\[-5pt]
&\multicolumn{3}{c@{}}{\hrulefill}\\
$\bolds{\sqrt{\frac{T}{\Delta_n}\var(\beta_{X,T}^{ts})}}$&\textbf{Median}&\textbf{IQR}&\textbf{MAD}\\
\hline
\multicolumn{4}{@{}l@{}}{\textit{Case A}}\\
\quad 3.3341&3.2774&0.2005&0.1005\\
\multicolumn{4}{@{}l@{}}{\textit{Case B}}\\
\quad 4.0320&3.8366&0.2638&0.1320\\
\multicolumn{4}{@{}l@{}}{\textit{Case C}}\\
\quad 4.9588&4.6678&0.2609&0.0590\\
\multicolumn{4}{@{}l@{}}{\textit{Case D}}\\
\quad 4.6626&4.7929&0.4596&0.2298\\
\hline
\end{tabular*}
\tabnotetext[]{tt1}{\textit{Notes}: $\var(\beta_{X,T}^{ts})$ is the exact variance of the
two-step estimator, computed from the $10{,}000$ Monte Carlo replications
of the estimator. $\widehat{\textrm{Ase}}(\beta_{X,T}^{ts}))$ is the
estimated asymptotic standard error using \refeq{eq:ts_sv_c}. MAD is
computed around the exact standard error of the estimator
$\sqrt{\frac{T}{\Delta_n}\var(\beta_{X,T}^{ts})}$.}
\end{table}

%s6 ###
\section{Proofs}\label{sec:proof}

The proof of Theorems~\ref{thms:cp_levy} and \ref{thms:cp_gen} follows
from results in \cite{TT07} and therefore is omitted here. For the rest
of the results, we first proof the ones for the \Lvb case, and then
proceed with those involving semimartingales with time-varying
characteristics. In what follows we use $\mathbb{E}^n_{i-1}$ and
$\mathbb{P}^n_{i-1}$ as a shorthand for
$\mathbb{E}(\cdot|\mathcal{F}_{(i-1)\Delta_n})$ and
$\mathbb{P}(\cdot|\mathcal{F}_{(i-1)\Delta_n})$, respectively. In the
proofs $K$ will denote a positive constant that does not depend on the
sampling frequency and might change from line to line.

%s6.1 ###
\subsection{\texorpdfstring{Proof of Theorem~\protect\ref{thms:asl}}{Proof of Theorem 3.3.}}

The proof of the theorem consists of showing (1) finite-dimensional
convergence (i.e., identifying the limit) and (2) tightness of the
sequence. In the proof we will show part (b) only. Part (a) can be
established in exactly the same way. We will assume that $A$ in
\refeq{eq:nu_1_2} is that of a standard stable process and therefore
$\Pi_{A,\beta}=1$. The result for an arbitrary $A$ then will follow
trivially by rescaling (and centering). In what follows $L$ will stand
for a standard symmetric $\beta$-stable process, defined on some
probability space which is possibly different from the original one.

\textit{Step \textup{1 (}Finite-dimensional convergence}). We start with establishing the final-dimensional
convergence. It will follow from Lemma~\ref{lema_fd} below in which we
denote with $\bullet$ the Hadamard product of two matrixes (i.e., the
element-by-element product). The stated lemma is slightly stronger than
what we need for two reasons. First, it contains locally uniform
convergence in $t$ and in the theorem we work with a fixed $T$. Second,
in the lemma we will show the finite-dimensional convergence for a
process $X$ defined in the following way:
%e6.1 ###
\begin{equation}\label{eq:js_00}
\qquad
X_t=\int_0^tm_{ds}\,ds+\int_0^t\int_{\mathbb{R}}\overline{\sigma}_{s-}\kappa(x)\tilde{\mu}(ds,dx)+\int_0^t\int_{\mathbb{R}}\overline{\sigma}_{s-}\kappa'(x)\mu(ds,dx),\hspace*{-10pt}
\end{equation}
where $\mu$ is the Poisson measure of Theorem~\ref{thms:asl};~for
arbitrary \cadlag processes $\sigma_{s}$ and $\widetilde{\sigma}_{s}$
with $K^{-1}<|\sigma_{s}|<K$ and $0\leq |\widetilde{\sigma}_{s}|\leq K$
for some $K>0$ and a Brownian motion $W_t$, $\overline{\sigma}_{s}$ is
defined via
$\overline{\sigma}_{s}=\sigma_{(i-1)\Delta_n}+\widetilde{\sigma}_{(i-1)\Delta_n}(W_s-W_{(i-1)\Delta_n})$
for $s\in [(i-1)\Delta_n,i\Delta_n)$ and further
$m_{ds}=m_{d,(i-1)\Delta_n}$ for $s\in [(i-1)\Delta_n,i\Delta_n)$.
Obviously $X_t$ includes the \Lvb case of Theorem~\ref{thms:asl}, and
the generalization will be needed later for the proof of
Theorem~\ref{thms:sv}.

\begin{lema}\label{lema_fd}
Let $\underline{\mathbf{p}}=(p_1,\ldots,p_k)'$ for some integer $k$,
$\mu_{\underline{\mathbf{p}}}=(\mu_{p_1},\ldots,\mu_{p_k})'$ and
$\mathbf{1}_k$ is $k\times 1$ vector of ones. Then, if $X$ is given by
\refeq{eq:js_00} and under the conditions of
Theorem~\textup{\ref{thms:asl}(b)}
(in particular all elements of $\underline{\mathbf{p}}$ are in
$[p_l,p_h]$), we have the following convergence locally uniformly in t:
%e6.2 ###
\begin{eqnarray}
\label{eq:js_0}&&\frac{1}{\sqrt{\Delta_n}}\widetilde{V}_t(\underline{\mathbf{p}},X,\Delta_n)\stackrel{\mathcal{L}-s}{\longrightarrow}\Xi(\underline{\mathbf{p}})_t,\nonumber
\\
\widetilde{V}_t(\underline{\mathbf{p}},X,\Delta_n)
&=&
\pmatrix{
\Delta_n^{\mathbf{1}_k-\underline{\mathbf{p}}/\beta}\bullet V_t(\underline{\mathbf{p}},X,2\Delta_n)-
\Delta_n^{\mathbf{1}_k-\underline{\mathbf{p}}/\beta}\bullet2^{\underline{\mathbf{p}}/\beta-\mathbf{1}_k}\cr
\hspace*{18pt}\bullet\mu_{\underline{\mathbf{p}}}(\beta)\bullet\displaystyle\sum_{i=1}^{[t/\Delta_n]}\biggl(\int_{(i-1)\Delta_n}^{i\Delta_n}|\overline{\sigma}_{s}|^{\beta}\,ds\biggr)^{\underline{\mathbf{p}}/\beta}\hspace*{-5pt}\vspace*{2pt}\cr
\hspace*{-8pt}\Delta_n^{\mathbf{1}_k-\underline{\mathbf{p}}/\beta}\bullet V_t(\underline{\mathbf{p}},X,\Delta_n)-\Delta_n^{\mathbf{1}_k-\underline{\mathbf{p}}/\beta}\bullet\mu_{\underline{\mathbf{p}}}(\beta)\cr
\hspace*{-6pt}\displaystyle\bullet\sum_{i=1}^{[t/\Delta_n]}\biggl(\int_{(i-1)\Delta_n}^{i\Delta_n}|\overline{\sigma}_{s}|^{\beta}\,ds\biggr)^{\underline{\mathbf{p}}/\beta}\hspace*{5pt}},
\\
V_t(\underline{\mathbf{p}},X,\iota\Delta_n)
&=&
(V_t(p_1,X,\iota\Delta_n),\ldots,V_t(p_k,X,\iota\Delta_n))',\qquad\iota=1,2,\nonumber
\end{eqnarray}
and the $\mathbb{R}^{2k}$-valued process
$\Xi(\underline{\mathbf{p}})_t$ is defined on an extension of the
original probability space, is continuous, and conditionally on the
$\sigma$-field $\mathcal{F}$ of the original probability space is
centered Gaussian with variance--covariance matrix process given by
$C_t$ defined via
%e6.3 ###
\begin{equation}\label{eq:js_3}
\qquad C_t(i,j)=
\cases{
\displaystyle\int_0^t|\sigma_{s}|^{p_i+p_j}\,ds\,2^{p_i/\beta+p_j/\beta-1}\bigl(\mu_{p_i+p_j}(\beta)-\mu_{p_i}(\beta)\mu_{p_j}(\beta)\bigr)\cr
\qquad\mbox{for $i=1,\ldots,k$; $j=1,\ldots,k$,}\cr
\displaystyle\int_0^t|\sigma_{s}|^{p_{i-k}+p_{j-k}}\,ds\,\bigl(\mu_{p_{i-k}+p_{j-k}}(\beta)-\mu_{p_{i-k}}(\beta)\mu_{p_{j-k}}(\beta)\bigr)\cr
\qquad\mbox{for $i=k+1,\ldots,2k$; $j=k+1,\ldots,2k$,}\cr
\displaystyle\int_0^t|\sigma_{s}|^{p_{i-k}+p_j}\,ds\,\bigl(\mu_{p_{i-k},p_j}(\beta)-2^{p_j/\beta}\mu_{p_{i-k}}(\beta)\mu_{p_j}(\beta)\bigr)\cr
\qquad\mbox{for $i=k+1,\ldots,2k$; $j=1,\ldots,k$}.}
\end{equation}
\end{lema}

\begin{pf}
We start with some notation. We set
$\widetilde{C}= C_t$ when $t=1$ and $\overline{\sigma}_{s}\equiv 1$ for
$\forall s\in[0,1]$. We further denote
%e6.4 ###
\begin{equation}\label{eq:js_000}
Y_t=\int_0^tm_{ds}\,ds+\int_0^t\!\!\int_{\mathbb{R}}\kappa(x)\tilde{\mu}(ds,dx)+\int_0^t\!\!\int_{\mathbb{R}}\kappa'(x)\mu(ds,dx),
\end{equation}
and
\begin{eqnarray}\label{eq:js_001}
X_t(\tau)
&=&
X_t-\sum_{s\leq t}\Delta X_s 1_{\{|\Delta X_s|<|\overline{\sigma}_{s-}|\tau\}},\nonumber
\\[-8pt]\\[-8pt]
Y_t(\tau)
&=&
Y_t-\sum_{s\leq t}\Delta Y_s 1_{\{|\Delta Y_s|<\tau\}},\qquad\tau>0.\nonumber
\end{eqnarray}
First, we have
\begin{eqnarray}
\Delta_n^{1-1/2-p_i/\beta}|V_t(p_i,X,\Delta_n)-V_t(p_i,X(\tau),\Delta_n)|
&\stackrel{\mathrm{u.c.p.}}{\longrightarrow}&
0,\qquad i=1,\ldots,k,\nonumber
\\[-8pt]\\[-8pt]
\qquad\hspace*{5pt}\Delta_n^{1-1/2-p_i/\beta}|V_t(p_i,X,2\Delta_n)-V_t(p_i,X(\tau),2\Delta_n)|
&\stackrel{\mathrm{u.c.p.}}{\longrightarrow}&
0,\qquad i=1,\ldots,k,\hspace*{-10pt}\nonumber
\end{eqnarray}
using the algebraic inequality $||a+b|^p-|a|^p|\leq |b|^p$ for $p\leq
1$ and the fact that $p_i<\beta/2$ for $i=1,\ldots,k$. Therefore we are
left with showing \refeq{eq:js_0} with
$V_t(\underline{\mathbf{p}},X,\Delta_n)$ and
$V_t(\underline{\mathbf{p}},X,2\Delta_n)$ substituted with
$V_t(\underline{\mathbf{p}},X(\tau),\Delta_n)$ and
$V_t(\underline{\mathbf{p}},X(\tau), 2\Delta_n)$, respectively.

For arbitrary power $p$ we set
\begin{eqnarray*}
\zeta(p)_i^n
&=&
(\zeta(p)_{i1}^n,\zeta(p)_{i2}^n)',\qquad i=1,2,\ldots,\biggl[\frac{t}{2\Delta_n}\biggr],
\\
\zeta(p)_{i1}^n
&=&
\Delta_n^{1/2}\biggl(\Delta_n^{-p/\beta}|\Delta_{2i-1}^nX(\tau)|^p+\Delta_n^{-p/\beta}|\Delta_{2i}^nX(\tau)|^p
\\
&&\hspace*{50pt}{}-
2\mu_p(\beta)\biggl(\frac{1}{\Delta_n}\int_{(i-1)\Delta_n}^{i\Delta_n}|\overline{\sigma}_{s}|^{\beta}\,ds\biggr)^{p/\beta}\biggr),
\\
\zeta(p)_{i2}^n
&=&
\Delta_n^{1/2}\biggl(\Delta_n^{-p/\beta}|\Delta_{2i-1}^nX(\tau)+\Delta_{2i}^nX(\tau)|^p
\\
&&\hphantom{\Delta_n^{1/2}\biggl(}{}-
2^{p/\beta}\mu_p(\beta)\biggl(\frac{1}{\Delta_n}\int_{(i-1)\Delta_n}^{i\Delta_n}|\overline{\sigma}_{s}|^{\beta}\,ds\biggr)^{p/\beta}\biggr).
\end{eqnarray*}
It is convenient also to write further
$\zeta(p)_{i1}^n=\xi(p)_{2i-1}+\xi(p)_{2i}$ with
\[
\xi(p)_j=\Delta_n^{1/2}\biggl(\Delta_n^{-p/\beta}|\Delta_{j}^nX(\tau)|^p
-
\mu_p(\beta)\biggl(\frac{1}{\Delta_n}\int_{(j-1)\Delta_n}^{j\Delta_n}|\overline{\sigma}_{s}|^{\beta}\,ds\biggr)^{p/\beta}\biggr)
\]
for $j=1,2,\ldots,2[\frac{t}{2\Delta_n}]$. Using Theorem IX.7.19 in
\cite{JS} it suffices to show the following for all $t>0$ and arbitrary
element $p$ from the vector $\underline{\mathbf{p}}$:
%e6.5 ###
\begin{equation}\label{eq:js_1}
\Biggl|\sum_{i=1}^{[t/(2\Delta_n)]}\mathbb{E}_{2i-2}^n(\zeta(p)_i^n)\Biggr|
\stackrel{\mathbb{P}}{\longrightarrow}
0,
\end{equation}
\begin{eqnarray}\label{eq:js_2}
&&
\sum_{i=1}^{[t/(2\Delta_n)]}\bigl(\mathbb{E}_{2i-2}^n[\zeta(p_q)_{is}^n\zeta(p_r)_{il}^n]-\mathbb{E}_{2i-2}^n(\zeta(p_q)_{is}^n)\mathbb{E}_{2i-2}^n(\zeta(p_r)_{il}^n)\bigr)\nonumber
\\[-8pt]\\[-8pt]
&&\qquad\stackrel{\mathbb{P}}{\longrightarrow}
C_t\bigl(q+(2-s)k,r+(2-l)k\bigr),\nonumber
\end{eqnarray}
where $s,l=1,2$ and $q,r=1,\ldots,k$,
%e6.6 ###
\begin{equation}\label{eq:js_4}
\sum_{i=1}^{[t/(2\Delta_n)]}\mathbb{E}_{2i-2}^n|\zeta(p)_i^n|^{2+\iota}\stackrel{\mathbb{P}}{\longrightarrow}0\qquad\mbox{for some }0<\iota<\beta/p-2,
\end{equation}
%e6.7 ###
\begin{equation}\label{eq:js_5}
\sum_{i=1}^{[t/(2\Delta_n)]}\mathbb{E}_{2i-2}^n[\zeta(p)_i^n(\Delta_{2i-1}^nM+\Delta_{2i}^nM)]\stackrel{\mathbb{P}}{\longrightarrow}0
\end{equation}
for $M$ being an arbitrary bounded local martingale defined on the
original probability space.

We start with \refeq{eq:js_1}. We prove it for the first element of
$\zeta(p)_i^n$ and arbitrary element $p$ of the vector
$\underline{\mathbf{p}}$, the proof for the second element of
$\zeta(p)_i^n$ is similar. Because of the assumption on the \Lvb
measure in \refeq{eq:nu} we can write
\begin{eqnarray*}
\mathbb{E}_{i-1}^n\biggl(|\Delta_n^{-1/\beta}\Delta_i^n X(\tau)|^p-\mu_p(\beta)
\biggl(\frac{1}{\Delta_n}\int_{(i-1)\Delta_n}^{i\Delta_n}|\overline{\sigma}_{s}|^{\beta}\,ds\biggr)^{p/\beta}\biggr)
=\sum_{j=1}^3A_{ij}^n
\end{eqnarray*}
for $i=1,2,\ldots,2[\frac{t}{2\Delta_n}]$ and where
\begin{eqnarray*}
A_{i1}^n
&=&
\mathbb{E}_{i-1}^n\biggl(\biggl|\Delta_n^{-1/\beta}\int_{(i-1)\Delta_n}^{i\Delta_n}\overline{\sigma}_{s-}\,dL_s\biggr|^p
-
\mu_p(\beta)\biggl(\frac{1}{\Delta_n}\int_{(i-1)\Delta_n}^{i\Delta_n}|\overline{\sigma}_{s}|^{\beta}\,ds\biggr)^{p/\beta}\biggr),
\\
A_{i2}^n
&=&
\mathbb{E}_{i-1}^n\biggl(\biggl|\Delta_n^{-1/\beta}\int_{(i-1)\Delta_n}^{i\Delta_n}\overline{\sigma}_{s-}\,dL_s+a_i\Delta_n^{-1/\beta}\biggr|^p
\\
&&\hspace*{67pt}{}-
\biggl|\Delta_n^{-1/\beta}\int_{(i-1)\Delta_n}^{i\Delta_n}\overline{\sigma}_{s-}\,dL_s\biggr|^p\biggr),
\\
A_{i3}^n
&=&
\mathbb{E}_{i-1}^n|\Delta_n^{-1/\beta}\Delta_i^n X(\tau)|^p
-
\mathbb{E}_{i-1}^n\biggl|\Delta_n^{-1/\beta}\int_{(i-1)\Delta_n}^{i\Delta_n}\overline{\sigma}_{s-}\,dL_s+a_i\Delta_n^{-1/\beta}\biggr|^p,
\end{eqnarray*}
with
%e6.8 ###
\begin{eqnarray} \label{eq:js_a}
a_i
&=&
m_{d,(i-1)\Delta_n}\Delta_n\nonumber
\\
&&{}-
\biggl(\int_{|x|>\tau}\kappa'(x)\nu_1(x)\,dx+2\int_{x:\nu_2(x)<0,|x|<\tau}\kappa(x)\nu_2(x)\,dx\biggr)
\\
&&\quad{}\times
\int_{(i-1)\Delta_n}^{i\Delta_n}\overline{\sigma}_{s}\,ds,\nonumber
\end{eqnarray}
where we recall that $L$ is a standard stable process which
is defined on an extension of the original probability space and is
independent of it. We have $a_i=0$ for $\beta\leq \sqrt{2}$, because of
our assumption of the symmetry of $\nu(x)$ and
$m_{d,(i-1)\Delta_n}\equiv 0$ for this case. Also, by the assumptions
of the theorem, $\beta'<\beta/2\leq 1$ and therefore the integral with
respect to $\nu_2$ in the definition of $a_i$ is well defined. Then,
using the algebraic inequality  $|x+y|^p\leq |x|^p+|y|^p$ for $p\leq 1$
and arbitrary $x$ and $y$, it is easy to show that for $A_{i3}^n$ we
have
\begin{eqnarray*}%\label{eq:an_00}
|A_{i3}^n|
&\leq&
K\Delta_n^{-p/\beta}\mathbb{E}_{i-1}^n\biggl|\int_{(i-1)\Delta_n}^{i\Delta_n}\overline{\sigma}_{s-}\,d\widetilde{L}^{(1)}_s\biggr|^p
+
K\Delta_n^{-p/\beta}\mathbb{E}_{i-1}^n\biggl|\int_{(i-1)\Delta_n}^{i\Delta_n}\overline{\sigma}_{s-}\,d\widetilde{L}^{(2)}_s\biggr|^p
\\
&&{}+
K\Delta_n^{-p/\beta}\mathbb{E}_{i-1}^n\biggl|\int_{(i-1)\Delta_n}^{i\Delta_n}\overline{\sigma}_{s-}\,d\widetilde{L}^{(3)}_s\biggr|^p,
\end{eqnarray*}
where $K$ is some constant and
%e6.9 ###
\begin{equation}\label{eq:L-s}
\qquad\cases{
\mbox{$\widetilde{L}^{(1)}$ is a pure-jump \Lvb process with \Lvb density of}\cr
\qquad\mbox{$-2\nu_2(x)1_{\{x:\nu_2(x)<0, |x|<\tau\}}$, zero drift and zero truncation function;}\cr
\mbox{$\widetilde{L}^{(2)}$ is a pure-jump \Lvb process with \Lvb density of}\cr
\qquad\mbox{$\nu_2(x)1_{\{x:\nu_2(x)>0, |x|<\tau\}}-\nu_2(x)1_{\{x:\nu_2(x)<0, |x|<\tau\}}$,}\cr
\qquad\mbox{zero drift and zero truncation function;}\cr
\mbox{$\widetilde{L}^{(3)}$ is a pure-jump \Lvb process with \Lvb density of}\cr
\qquad\mbox{$\nu_1(x)1_{\{|x|>\tau\}}$, zero drift and zero truncation function.}
}\hspace*{-10pt}
\end{equation}
The three processes are well defined because $\beta'<1$ and are defined
on an extension of the original probability space and independent from
the original filtration. Then, using the fact that
$\overline{\sigma}_{s-}$ is independent from the processes
$\widetilde{L}^{(i)}$ for $i=1,2,3$,
$\mathbb{E}|\overline{\sigma}_{s}|^p<\infty$ for
$s\in[(i-1)\Delta_n,i\Delta_n)$ and any positive $p$, the H\"{o}lder's
inequality, and the basic one $|\sum_{i}|a_i||^p\leq \sum_{i}|a_i|^p$
for $p\leq 1$ and arbitrary $a_i$, we easily have
%e6.10 ###
\begin{equation}\label{eq:an_1}
|A_{i3}^n|\leq K\Delta_n^{p/\beta'\wedge 1-p/\beta-\iota}
\end{equation}
for any $\iota>0$. Taking into account the restriction on $p$ and
$\beta'$, we have $p/\beta'\wedge 1-p/\beta-\iota>1/2$ for some
$\iota>0$. In a similar way we can show $|\widetilde{A}_{i3}^n|\leq
K\Delta_n^{1/2+\iota}$ for some $\iota>0$ where
\begin{eqnarray*}
\widetilde{A}_{i3}^n
&=&
\mathbb{E}_{i-1}^n\biggl(|\Delta_n^{-1/\beta}\Delta_i^n X(\tau)|^p\Delta_i^nW
\\
&&\hphantom{\mathbb{E}_{i-1}^n\biggl(}{}-
\biggl|\Delta_n^{-1/\beta}\int_{(i-1)\Delta_n}^{i\Delta_n}\overline{\sigma}_{s-}\,dL_s+a_i\Delta_n^{-1/\beta}\biggr|^p\Delta_i^nW\biggr).
\end{eqnarray*}
Further, since
$\int_{(i-1)\Delta_n}^{i\Delta_n}\overline{\sigma}_{s-}\,dL_s
\stackrel{d}{=}L_{b_{i,n}}$ for
$b_{i,n}=\int_{(i-1)\Delta_n}^{i\Delta_n}|\overline{\sigma}_{s}|^{\beta}\,ds$,
and using the self-similarity property of a strictly stable process, we
have $A_{i1}^n=0$. We have similarly $\widetilde{A}_{i1}^n=0$, where
\begin{eqnarray*}
\widetilde{A}_{i1}^n
&=&
\mathbb{E}_{i-1}^n\biggl(\biggl|\Delta_n^{-1/\beta}\int_{(i-1)\Delta_n}^{i\Delta_n}\overline{\sigma}_{s-}\,dL_s\biggr|^p\Delta_i^nW
\\
&&\hphantom{\mathbb{E}_{i-1}^n\biggl(}{}-
\mu_p(\beta)\biggl(\frac{1}{\Delta_n}\int_{(i-1)\Delta_n}^{i\Delta_n}|\overline{\sigma}_{s}|^{\beta}\,ds\biggr)^{p/\beta}\Delta_i^nW\biggr),
\end{eqnarray*}
because $W$ is independent from $L$. Next, to prove
\refeq{eq:js_1}, we need only show that $|A_{i2}^n|\leq
K\Delta_n^{1/2+\iota}$ for some $\iota>0$.  We show this only for the
case $\beta>\sqrt{2}$, since for $\beta\leq\sqrt{2}$ it is trivially
satisfied. For the proof we make use of the following general
inequality for arbitrary real numbers $x$ and $y$ and $p\leq 1$:
\begin{eqnarray}\label{eq:an_0}
&&
\bigl||x+y|^p-|x|^p-p|x|^{p-1}\operatorname{sign}\{x\}y1_{\{|x|\neq 0,|y|\leq|x|/2\}}\bigr|\nonumber
\\[-8pt]\\[-8pt]
&&\qquad{}\leq
K\frac{|y|^{p+1-\iota}}{|x|^{1-\iota}}1_{\{|x|\neq 0\}}+|y|^p1_{\{|x|= 0\cup|y|>|x|/2\}}\nonumber
\end{eqnarray}
for some $\iota>0$ and a positive constant $K$. The inequality follows
by looking at the difference $|x+y|^p-|x|^p$ on the following two sets:
$|y|\leq |x|/2$ and $|y|>|x|/2$. On the former we apply a second-order
Taylor series approximation and further use $|y|/|x|\leq 1/2$ on this
set [therefore \refeq{eq:an_0} holds with $K=2^{p-2-\iota}p(1-p)$]. On
the set $|y|>|x|/2$ we use the subadditivity of the function $|x|^p$.
We can substitute in the above inequality $x$ with
$\Delta_n^{-1/\beta}L_{b_{i,n}}$ and $y$ with $a_i\Delta_n^{-1/\beta}$.
Then, by first conditioning on the filtration generated by
$\overline{\sigma}_{s}$, and then using the fact that $L$ has symmetric
distribution, we get
%e6.11 ###
\begin{equation}\label{eq:an_2}
\qquad\mathbb{E}_{i-1}^n\bigl(|\Delta_n^{-1/\beta}L_{b_{i,n}}|^{p-1}\operatorname{sign}\{L_{b_{i,n}}\}a_i\Delta_n^{-1/\beta}1_{\{|L_{b_{i,n}}|\neq 0,|L_{b_{i,n}}|\geq 2|a_i|\}}\bigr)=0.
\end{equation}
Next we have for some $p_0,p_1>0$ (note that we have universal bounds
on $\sigma_{s}$ and $\widetilde{\sigma}_{s}$)
%e6.12 ###
\begin{equation}\label{eq:an_2a}
\mathbb{E}_{i-1}^n\biggl(\int_{(i-1)\Delta_n}^{i\Delta_n}|\overline{\sigma}_{s}|^{p_0}\,ds\biggr)^{-p_1}
\leq
K\mathbb{E}_{i-1}^n(T_b\wedge \Delta_n)^{-p_1}
<
K\Delta_n^{-p_1},
\end{equation}
where $T_b$ is the hitting time of  the Brownian motion
$(W_s-W_{(i-1)\Delta_n})_{s\geq (i-1)\Delta_n}$ of the level $b$ for
$b=-\sigma_{(i-1)\Delta_n}/(2K)\neq 0$ for some positive $K$, whose
negative powers (of $T_b$) are finite. Then for $\iota$ such that
$0<\iota<p-\frac{2-\beta}{2(\beta-1)}$ (recall the assumption on $p$
for $\beta>\sqrt{2}$) we have
%e6.13 ###
\begin{eqnarray}\label{eq:an_3}
&&
\mathbb{E}_{i-1}^n\biggl(\frac{|a_i\Delta_n^{-1/\beta}|^{p+1-\iota}}{|\Delta_n^{-1/\beta}L_{b_{i,n}}|^{1-\iota}}1_{\{L_{b_{i,n}}\neq 0\}}\biggr)\nonumber
\\
&&\qquad\leq
\mathbb{E}_{i-1}^n[|a_i\Delta_n^{-1/\beta}|^{p+1-\iota}|\Delta_n^{-1/\beta}b_{i,n}^{1/\beta}|^{\iota-1}]\mathbb{E}(|L_1|^{\iota-1})
\\
&&\qquad\leq
K \Delta_n^{1/2+\iota'},\nonumber
\end{eqnarray}
with some $\iota'>0$ and a positive constant $K$. This follows from the
self-similarity of the strictly stable process, the fact that
$\mathbb{E}|L_1|^{1-\iota}$ $<$ $\infty$ since $\iota\in(0,1)$ (see,
e.g., \cite{SATO}) and the preceding inequality \refeq{eq:an_2a}.
Similarly, for some $\iota\in (0,p-\frac{2-\beta}{2(\beta-1)})$ using
the Chebyshev's inequality we have
\begin{eqnarray}\label{eq:an_3a}
\mathbb{E}_{i-1}^n|a_i\Delta_n^{-1/\beta}|^p1_{\{|L_{b_{i,n}}|<2|a_i|\}}
&\leq&
K\mathbb{E}(|L_1|^{\iota-1})\Delta_n^{(1-1/\beta)(p+1-\iota)}\nonumber
\\[-8pt]\\[-8pt]
&\leq&
K\Delta_n^{1/2+\iota'}\nonumber
\end{eqnarray}
with some $\iota'>0$. Combining \refeq{eq:an_1}--\refeq{eq:an_3a} and
using that stable distribution has a density with respect to Lebesgue
measure (see, e.g., Remark 14.18 in \cite{SATO}) we prove
$|A_{i2}^n|\leq K\Delta_n^{1/2+\iota}$ for some $\iota>0$ and thus
\refeq{eq:js_1} follows. Similarly we have $|\widetilde{A}_{i2}^n|\leq
K\Delta_n^{1/2+\iota}$ for some $\iota>0$ where
\begin{eqnarray*}
\widetilde{A}_{i2}^n
&=&
\mathbb{E}_{i-1}^n\biggl(\biggr|\Delta_n^{-1/\beta}\int_{(i-1)\Delta_n}^{i\Delta_n}\overline{\sigma}_{s-}\,dL_s+a_i\Delta_n^{-1/\beta}\biggr|^p\Delta_i^n W
\\
&&\hspace*{67pt}{}-
\biggl|\Delta_n^{-1/\beta}\int_{(i-1)\Delta_n}^{i\Delta_n}\overline{\sigma}_{s-}\,dL_s\biggr|^p\Delta_i^n W\biggr).
\end{eqnarray*}
%%%%%%%%%%%%%%%%%%%%%%%%%%%%%%%%%%%%%%%%%%%%%%%%%%%%%%%%%%%%%%%%%%%%%%%%%%%%%%%%%%%%%%%%%%%%%%%%%%%%%%%%%%%%%%%%%%%%%%%%%%%%%%%%%%%%%%%%%%%%%%%%%%%%%%%%%%%%%%%%

Before proceeding with \refeq{eq:js_2} we derive a result that we make
use of later for the proof of Theorem~\ref{thms:sv}. First, for two
random variables $X_1$ and $X_2$ and some $\varepsilon>0$ we have
%e6.14 ###
\begin{equation}
\mathbb{P}(|X_1+X_2|\leq \varepsilon)\leq \mathbb{P}(|X_1|\geq\varepsilon)+\mathbb{P}(|X_2|\leq 2\varepsilon).
\end{equation}
Then we can apply this inequality twice, use the fact that
$\int_{[-1,1]}|x|^{\beta'+\alpha'}\nu_2(x)\,dx$ $<$ $\infty$ for any
$\alpha'>0$, the fact that $|\Delta X_s(\tau)|\leq
\tau|\overline{\sigma}_{s-}|$; the fact that the stable distribution
has finite moments for powers that are negative but higher than $-1$;
the bound in \refeq{eq:an_2a} and finally the Chebyshev's inequality
to get
%e6.15 ###
\begin{eqnarray}\label{eq:an_3b}
\mathbb{P}_{i-1}^n\bigl(\Delta_n^{-1/\beta}|\Delta_i^n X(\tau)|\leq\varepsilon\bigr)
&\leq&
\sum_{j=1}^3\mathbb{P}_{i-1}^n\biggl(\biggl|\int_{(i-1)\Delta_n}^{i\Delta_n}\overline{\sigma}_{s-}\,d\widetilde{L}^{(j)}_s\biggr|\geq0.5\Delta_n^{1/\beta}\varepsilon\biggr)\nonumber
\\
&&{}+
\mathbb{P}_{i-1}^n(|a_i\Delta_n^{-1/\beta}+\Delta_n^{-1/\beta}L_{b_{i,n}}|\leq4\varepsilon)
\\
&\leq&
K\biggl(\varepsilon^{\alpha}+\Delta_n^{(1-1/\beta)\alpha}+\frac{\Delta_n^{p/\beta'-p/\beta-\alpha'}}{\varepsilon^{p}}\biggr)\nonumber
\end{eqnarray}
for any $\alpha\in (0,1)$, $p\leq\beta'$ and $\alpha'>0$ and where $K$
is some positive constant that \textit{does not depend} on $\varepsilon$.
Similarly for two random variables $X_1$ and $X_2$ and $p>0$ and
$\varepsilon>0$ we can derive
\begin{eqnarray*}%\label{eq:an_3c}
&&
\mathbb{E}\bigl(|X_1+X_2|^{-p}1_{\{|X_1+X_2|\geq \varepsilon\}}\bigr)
\\
&&\qquad\leq
K\bigl[\varepsilon^{-p}\mathbb{P}(|X_2|\geq k\varepsilon)
+
\mathbb{E}\bigl(|X_1|^{-p}1_{\{|X_1|>(1-k)\varepsilon\}}\bigr)\bigr]
\end{eqnarray*}
for any $k\in (0,1)$ and where the constant $K$ depends on $k$ only.
Using this inequality then it is easy to derive the following bound:
\begin{eqnarray}\label{eq:an_3d}
&&
\mathbb{E}_{i-1}^n\bigl(|\Delta_n^{-1/\beta}\Delta_i^n X(\tau)|^{-p}1_{\{\Delta_n^{-1/\beta}|\Delta_i^n X(\tau)|\geq\varepsilon\}}\bigr)\nonumber
\\[-8pt]\\[-8pt]
&&\qquad\leq
K\biggl(\varepsilon^{(1-p)\wedge 0-\alpha'}+\frac{\Delta_n^{1-\beta'/\beta-\alpha'}}{\varepsilon^{p+\beta'}}\biggr)\nonumber
\end{eqnarray}
for any $p,\alpha'>0$ and where the constant $K$ does not depend on
$\varepsilon$.
%%%%%%%%%%%%%%%%%%%%%%%%%%%%%%%%%%%%%%%%%%%%%%%%%%%%%%%%%%%%%%%%%%%%%%%%%%%%%%%%%%%%%%%%%%%%%%%%%%%%%%%%%%%%%%%%%%%%%%%%%%%%%%%%%%%%%%%%%%%%%%%%%%%%%%%%%%%%%%%%

We continue with \refeq{eq:js_2}. First using Lemma 1(b) in
\cite{TT07}, since for each element $p$ of the vector
$\underline{\mathbf{p}}$ we have $2p<\beta$, we have [recall the
notation in \refeq{eq:js_000} and \refeq{eq:js_001}]
\begin{eqnarray*}%\label{eq:an_4a}
&&
\mathbb{E}_{i-1}^n|\Delta_n^{-1/\beta}\Delta_i^nY(\tau)|^{p_q+p_r}-\mathbb{E}_{i-1}^n|\Delta_n^{-1/\beta}\Delta_i^nY(\tau)|^{p_q}
\mathbb{E}_{i-1}^n|\Delta_n^{-1/\beta}\Delta_i^nY(\tau)|^{p_r}
\\
&&\qquad\stackrel{\mathbb{P}}{\longrightarrow}
\widetilde{C}(k+q,k+r),
\end{eqnarray*}
\begin{eqnarray*}%\label{eq:an_4b}
&&
\tfrac{1}{2}\mathbb{E}_{2i-2}^n|\Delta_n^{-1/\beta}\Delta_{2i-1}^nY(\tau)+\Delta_n^{-1/\beta}\Delta_{2i}^nY(\tau)|^{p_q+p_r}
\\
&&\quad{}-
\tfrac{1}{2}\mathbb{E}_{2i-2}^n|\Delta_n^{-1/\beta}\Delta_{2i-1}^nY(\tau)+\Delta_n^{-1/\beta}\Delta_{2i}^nY(\tau)|^{p_q}
\\
&&\quad\hphantom{-}{}\times
\mathbb{E}_{2i-2}^n|\Delta_n^{-1/\beta}\Delta_{2i-1}^nY(\tau)+\Delta_n^{-1/\beta}\Delta_{2i}^nY(\tau)|^{p_r}
\\
&&\qquad\stackrel{\mathbb{P}}{\longrightarrow}
\widetilde{C}(q,r),
\end{eqnarray*}
\begin{eqnarray*}%\label{eq:an_4c}
&&
\mathbb{E}_{2i-2}^n|\Delta_n^{-1/\beta}\Delta_{2i-1}^nY(\tau)+\Delta_n^{-1/\beta}\Delta_{2i}^nY(\tau)|^{p_q}|\Delta_n^{-1/\beta}\Delta_{2i-1}^nY(\tau)|^{p_r}
\\
&&\quad{}-
\mathbb{E}_{2i-2}^n|\Delta_n^{-1/\beta}\Delta_{2i-1}^nY(\tau)+\Delta_n^{-1/\beta}\Delta_{2i}^nY(\tau)|^{p_q}\mathbb{E}_{i-1}^n|\Delta_n^{-1/\beta}\Delta_i^nY(\tau)|^{p_r}
\\
&&\qquad\stackrel{\mathbb{P}}{\longrightarrow}\widetilde{C}(q,k+r),
\end{eqnarray*}
where  $q,r=1,\ldots,k$ and for the first limit
$i=1,2,\ldots,2[\frac{t}{2\Delta_n}]$ while for the last two
$i=1,2,\ldots,[\frac{t}{2\Delta_n}]$. Next, by Riemann integrability, we
have
%e6.16 ###
\begin{equation}\label{eq:an_4d}
\Delta_n\sum_{i=1}^{[t/\Delta_n]}\bigl|\sigma_{(i-1)\Delta_n}\bigr|^p
\stackrel{\mathbb{P}}{\longrightarrow}
\int_0^t|\sigma_{s}|^p\,ds,\qquad p>0.
\end{equation}
Therefore, to show \refeq{eq:js_3} we need only to prove that for
arbitrary $p<\beta$
%e6.17 ###
\begin{equation}\label{eq:an_4e}
\mathbb{E}_{i-1}^n\bigl||\Delta_n^{-1/\beta}X(\tau)|^p-\bigl|\Delta_n^{-1/\beta}\sigma_{(i-1)\Delta_n}Y(\tau)\bigr|^p\bigr|\leq K\Delta_n^{\iota}
\end{equation}
for some $\iota>0$. But this follows by using the
Burkholder--Davis--Gundy inequality (if $\beta>1$) and the elementary
one $(\sum_i|a_i|)^p\leq \sum_i|a_i|^p$ for arbitrary reals $a_i$ and
some $p\leq 1$, together with the definition of the process
$\overline{\sigma}_{s}$.

Turning to \refeq{eq:js_4}, we show it only for the first component of
$\zeta(p)_i^n$, the proof for the second one being exactly the same.
Using again Lemma 1(b) in \cite{TT07} we have
\begin{eqnarray}
&&
\mathbb{E}_{i-1}^n\bigl(\Delta_n^{-(2+\iota)p/\beta}|\Delta_i^nY(\tau)|^{(2+\iota)p}\bigr)\nonumber
\\[-8pt]\\[-8pt]
&&\qquad\stackrel{\mathbb{P}}{\longrightarrow}
\mathbb{E}\bigl(|L_1|^{(2+\iota)p}\bigr)\nonumber
\end{eqnarray}
for $i=1,2,\ldots,2[\frac{t}{2\Delta_n}]$ and $0<\iota<\beta/p-2$. Then
\refeq{eq:js_4} follows by combining this result with
\refeq{eq:an_4d}--\refeq{eq:an_4e}. We are left with proving
\refeq{eq:js_5}. It suffices to show
\begin{eqnarray}\label{eq:js_6}
&&
\qquad\Delta_n^{1/2}\sum_{i=1}^{2[t/(2\Delta_n)]}\mathbb{E}_{i-1}^n\biggl(\Delta_n^{-p/\beta}|\Delta_i^nX(\tau)|^p\Delta_i^nM\nonumber
\\[-8pt]\\[-8pt]
&&\qquad\hphantom{\Delta_n^{1/2}\sum_{i=1}^{2[t/(2\Delta_n)]}\mathbb{E}_{i-1}^n\biggl(}{}-
\mu_p(\beta)\biggl(\frac{1}{\Delta_n}\int_{(i-1)\Delta_n}^{i\Delta_n}|\overline{\sigma}_{s}|^{\beta}\,ds\biggr)^{p/\beta}\Delta_i^nM\biggr)
\stackrel{\mathbb{P}}{\longrightarrow}
0.\nonumber
\end{eqnarray}
First, if $M$ is a discontinuous martingale, then using
\refeq{eq:js_1}--\refeq{eq:js_4}, we have that
$\sum_{i=1}^{2[t/(2\Delta_n)]}\xi(p)_i^n$ is C-tight, that is, it is
tight and any limit is continuous. At the same time
$\sum_{i=1}^{2[t/(2\Delta_n)]}\Delta_i^n M$ trivially converges to a
discontinuous limit. Therefore the pair
$(\sum_{i=1}^{2[t/(2\Delta_n)]}\xi(p)_i^n,\sum_{i=1}^{2[t/(2\Delta_n)]}\Delta_i^n
M)$ is tight (see \cite{JS}, Theorem VI3.33(b)). But then the left-hand
side of \refeq{eq:js_6} converges to the predictable version of the
quadratic covariation of the limits of
$\sum_{i=1}^{2[t/(2\Delta_n)]}\xi(p)_i^n$ and
$\sum_{i=1}^{2[t/(2\Delta_n)]}\Delta_i^n M$ (use Theorem VI.6.29 of \cite{JS} for this), which is zero since
continuous and discontinuous martingales are orthogonal (see \cite{JS},
Proposition I.4.15).

Second if $M$ is a continuous martingale orthogonal to the Brownian
motion $W_t$ used in defining $\overline{\sigma}_{t}$, we can proceed
similarly to \cite{BNGJPS05} and argue as follows. If we set
$N_t=\mathbb{E}(|\Delta_i^nX(\tau)|^p|\mathcal{F}_t)$ for $t\geq
(i-1)\Delta_n$, then $(N_t)_{t\geq (i-1)\Delta_n}$ is a martingale. It
remains also martingale, conditionally on
$\mathcal{F}_{(i-1)\Delta_n}$, for the filtration generated by the
Poisson measure $\mu$ and the Brownian motion
$(W_t-W_{(i-1)\Delta_n})_{t\geq (i-1)\Delta_n}$ since $\Delta_i^nX$ is
uniquely determined by these processes. Therefore, by a martingale
representation theorem (see \cite{JS}, Theorem III.4.34)
\begin{eqnarray*}
N_t
&=&
N_{(i-1)\Delta_n}+\int_{(i-1)\Delta_n}^{t}\int_{\mathbb{R}}\delta'(s,x)\widetilde{\mu}(ds,dx)
\\
&&{}+
\int_{(i-1)\Delta_n}^{t}\eta_s\,dW_s,
\end{eqnarray*}
when $t\geq(i-1)\Delta_n$ for an appropriate predictable function
$\delta'(s,x)$ and process $\eta_s$. Therefore $N_t$ is a sum of
pure-discontinuous martingale, which hence is orthogonal to
$M_t-M_{(i-1)\Delta_n}$ (see \cite{JS}, Definition I.4.11), and a continuous
martingale which is also orthogonal to $M_t-M_{(i-1)\Delta_n}$ because
of our assumption on $M$.  This implies that for $M$ a continuous
martingale orthogonal to the Brownian motion we have
\begin{eqnarray*}
&&
\mathbb{E}_{i-1}^n\biggl(\biggl[\Delta_n^{-p/\beta}|\Delta_i^nX(\tau)|^p
-
\mu_p(\beta)\biggl(\frac{1}{\Delta_n}\int_{(i-1)\Delta_n}^{i\Delta_n}|\overline{\sigma}_{s}|^{\beta}\,ds\biggr)^{p/\beta}\biggr]\Delta_i^nM\biggr)
\\
&&\qquad=
\mathbb{E}_{i-1}^n(\Delta_i^nN\Delta_i^nM)=0,
\end{eqnarray*}
and this shows \refeq{eq:js_6} in this case.

The only case that remains to be covered is when $M=W$. For this case
we can use the bounds derived above for $\widetilde{A}_{i1}$,
$\widetilde{A}_{i2}$ and $\widetilde{A}_{i3}$ and from here
\refeq{eq:js_6} follows easily in this case.
\end{pf}

\textit{Step \textup{2} (Tightness)}. We are left with
establishing tightness, which follows from the next lemma.

\begin{lema}\label{lema_ti}
Assume that $X$ is given by \refeq{eq:js_00} and that the conditions of
Theorem~\ref{thms:asl} hold. Then for a \textit{fixed} $T>0$ we have
that the sequence
\[
\frac{1}{\sqrt{\Delta_n}}\widetilde{V}_T(\underline{\mathbf{p}},X,\Delta_n)
\]
for $\widetilde{V}_T(\underline{\mathbf{p}},X,\Delta_n)$ defined in
\refeq{eq:js_0}, is tight on the space of continuous functions
$\mathcal C(\mathbb{R}^2,[p_l,p_h])$ equipped with the uniform
topology, where $p_l$ and $p_h$ satisfy the conditions of part \textup{(b)} of
Theorem~\ref{thms:asl}.
\end{lema}

\begin{pf}
We will prove only that the sequence
\begin{eqnarray*}
\widehat{V}_T(p,X,\Delta_n)
&=&
\Delta_n^{1/2-p/\beta}V_T(p,X,\Delta_n)
\\
&&{}-
\Delta_n^{1/2}\mu_{p}(\beta)\sum_{i=1}^{[T/\Delta_n]}\biggl(\frac{1}{\Delta_n}\int_{(i-1)\Delta_n}^{i\Delta_n}|\overline{\sigma}_{s}|^{\beta}\,ds\biggr)^{p/\beta}
\end{eqnarray*}
is tight in the space of $\mathbb{R}$-valued functions on $[p_l,p_h]$
and the arguments generalize to the tightness of
$\frac{1}{\sqrt{\Delta_n}}\widetilde{V}_T(\underline{\mathbf{p}},X,\Delta_n)$.
For arbitrary $p_l\leq p<q \leq p_h$ we can write
\[
|\widehat{V}_T(q,X,\Delta_n)-\widehat{V}_T(p,X,\Delta_n)|\leq\sum_{i=1}^4A_i^n(p,q),
\]
where
\begin{eqnarray*}
A_1^n(p,q)
&=&
\Delta_n^{-1/2}\bigl|\Delta_n^{1-q/\beta}\bigl(V_T(q,X,\Delta_n)-V_T(q,X(\tau),\Delta_n)\bigr)
\\
&&\hphantom{\Delta_n^{-1/2}\bigl|}{}-
\Delta_n^{1-p/\beta}\bigl(V_T(p,X,\Delta_n)-V_T(p,X(\tau),\Delta_n)\bigr)\bigr|,
\end{eqnarray*}
and for $i=2,3,4$, $A_i^n(p,q)\stackrel{d}{=}\widetilde{A}_i^n(p,q)$
with
\begin{eqnarray*}
\widetilde{A}_2^n(p,q)
&=&
\Delta_n^{1/2}\Biggl|\sum_{i=1}^{[T/\Delta_n]}\biggl[\biggl|\Delta_n^{-1/\beta}\int_{(i-1)\Delta_n}^{i\Delta_n}
\overline{\sigma}_{s-}\,dL_s\biggr|^q
\\
&&\hphantom{\Delta_n^{1/2}\Biggl|\sum_{i=1}^{[T/\Delta_n]}\biggl[}{}-
\mu_q(\beta)\biggl(\frac{1}{\Delta_n}\int_{(i-1)\Delta_n}^{i\Delta_n}|\overline{\sigma}_{s}|^{\beta}\,ds\biggr)^{q/\beta}
\\
&&\hphantom{\Delta_n^{1/2}\Biggl|\sum_{i=1}^{[T/\Delta_n]}\biggl[}{}-
\biggl|\Delta_n^{-1/\beta}\int_{(i-1)\Delta_n}^{i\Delta_n}\overline{\sigma}_{s-}\,dL_s\biggr|^p
\\
&&\hphantom{\Delta_n^{1/2}\Biggl|\sum_{i=1}^{[T/\Delta_n]}\biggl[}{}+
\mu_p(\beta)\biggl(\frac{1}{\Delta_n}\int_{(i-1)\Delta_n}^{i\Delta_n}|\overline{\sigma}_{s}|^{\beta}\,ds\biggr)^{p/\beta}\biggr]\Biggr|,
\\
\widetilde{A}_3^n(p,q)
&=&
\Delta_n^{1/2}\Biggl|\sum_{i=1}^{[T/\Delta_n]}\biggl[\biggl|\Delta_n^{-1/\beta}\int_{(i-1)\Delta_n}^{i\Delta_n}\overline{\sigma}_{s-}\,dL_s+a_i\Delta_n^{-1/\beta}\biggr|^q
\\
&&\hphantom{\Delta_n^{1/2}\Biggl|\sum_{i=1}^{[T/\Delta_n]}\biggl[}{}-
\biggl|\Delta_n^{-1/\beta}\int_{(i-1)\Delta_n}^{i\Delta_n}\overline{\sigma}_{s-}\,dL_s\biggr|^q
\\
&&\hphantom{\Delta_n^{1/2}\Biggl|\sum_{i=1}^{[T/\Delta_n]}\biggl[}{}-
\biggl|\Delta_n^{-1/\beta}\int_{(i-1)\Delta_n}^{i\Delta_n}\overline{\sigma}_{s-}\,dL_s+a_i\Delta_n^{-1/\beta}\biggr|^p
\\
&&\hspace*{107pt}{}+
\biggl|\Delta_n^{-1/\beta}\int_{(i-1)\Delta_n}^{i\Delta_n}\overline{\sigma}_{s-}\,dL_s\biggr|^p\biggr]\Biggr|,
\end{eqnarray*}
where $a_i$ is defined in \refeq{eq:js_a} in the proof of
Lemma~\ref{lema_fd} and $\widetilde{A}_4^n(p,q)$ is a residual term
whose moments involve the processes $\widetilde{L}^{(1)}$,
$\widetilde{L}^{(2)}$ and $\widetilde{L}^{(3)}$ of \refeq{eq:L-s}. It
can be shown using the continuity of the power function and the
restriction on $\nu_2(x)$ that
%e6.18 ###
\begin{eqnarray}
\limsup_{\Delta_n\downarrow 0}\mathbb{E}\Bigl(\sup_{p,q\in[p_l,p_h]}\widetilde{A}_4^n(p,q)\Bigr)=0.
\end{eqnarray}
For $A_1^n(p,q)$ we can first apply the inequality $||a+b|^p-|a|^p|\leq
|b|^p$ for $p\leq 1$, and then use the continuity of the power function
for positive powers to show that
%e6.19 ###
\begin{equation}\label{eq:at0}
\sup_{p,q\in[p_l,p_h]} A_1^n(p,q)\stackrel{\mathrm{a.s.}}{\longrightarrow}0.
\end{equation}
For $\widetilde{A}_2^n(p,q)$ we easily have for $p,q \in[p_l,p_h]$
%e6.20 ###
\begin{equation}\label{eq:at1}
\mathbb{E}(\widetilde{A}_2^n(p,q))^2\leq K(p-q)^2,
\end{equation}
and Theorem 12.3 in \cite{Billingsley} implies tightness.
%of $A_2^n(p,q)$
Turning to  $\widetilde{A}_3^n(p,q)$, it is identically $0$ for
$\beta\leq \sqrt{2}$ due to our assumptions. So we look at the case
$\beta>\sqrt{2}$. We can decompose $\widetilde{A}_3^n(p,q)$ as
$\widetilde{A}_3^n(p,q)\leq
\widetilde{A}_{31}^n(p,q)+\widetilde{A}_{32}^n(p,q)$ with
\[
\cases{
\displaystyle\widetilde{A}_{31}^n(p,q)=\Delta_n^{1/2}\Biggl|\sum_{i=1}^{[T/\Delta_n]}[c_i(q)-c_i(p)]1_{\{C_i^n\}}\Biggr|,\cr
\displaystyle\widetilde{A}_{32}^n(p,q)=\Delta_n^{1/2}\Biggl|\sum_{i=1}^{[T/\Delta_n]}[c_i(q)-c_i(p)]1_{\{(C_i^{n})^c\}}\Biggr|,
}
\]
where $C_i^n=\{|\int_{(i-1)\Delta_n}^{i\Delta_n}
\overline{\sigma}_{s-}\,dL_s|\neq 0, |\int_{(i-1)\Delta_n}^{i\Delta_n}
\overline{\sigma}_{s-}\,dL_s|\geq 2|a_i|\}$ and
\[
c_i(p)=\biggl|\Delta_n^{-1/\beta}\int_{(i-1)\Delta_n}^{i\Delta_n}\overline{\sigma}_{s-}\,dL_s+a_i\Delta_n^{-1/\beta}\biggr|^p
-
\biggl|\Delta_n^{-1/\beta}\int_{(i-1)\Delta_n}^{i\Delta_n}\overline{\sigma}_{s-}\,dL_s\biggr|^p.
\]

For $\widetilde{A}_{31}^n(p,q)$ we can write
\begin{eqnarray}\label{eq:at2}
\mathbb{E}(\widetilde{A}_{31}^n(p,q))^2
&\leq&
K\mathbb{E}\bigl([c_i(q)-c_i(p)]^21_{\{C_i^n\}}\bigr)\nonumber
\\[-8pt]\\[-8pt]
&&{}+
K\Delta_n\Biggl(\sum_{i=1}^{[T/\Delta_n]}\mathbb{E}_{i-1}^n\bigl([c_i(q)-c_i(p)]1_{\{C_i^n\}}\bigr)\Biggr)^2.\nonumber
\end{eqnarray}
For the first expectation on the left-hand side of \refeq{eq:at2} we
have, similarly to \refeq{eq:at1},
%e6.21 ###
\begin{equation}\label{eq:at3}
\mathbb{E}\bigl([c_i(q)-c_i(p)]^21_{\{C_i^n\}}\bigr)\leq K(p-q)^2.
\end{equation}
For the second expectation on the right-hand side of \refeq{eq:at2}, we
apply the following inequality, similarly to \refeq{eq:an_0}. For every
$x$ and $y$ and $p,q\in[p_l; p_h]$ we have
\begin{eqnarray*}%\label{eq:at4}
&&
\bigl||x+y|^p-|x|^p-|x+y|^q+|x|^q
\\
&&\quad{}-
(p|x|^{p-1}-q|x|^{q-1})\operatorname{sign}\{x\}y1_{\{|x|\neq 0, 2|y|\leq |x|\}}\bigr|1_{\{|x|\neq 0, 2|y|\leq |x|\}}
\\
&&\qquad\leq
K|p-q|\frac{(|y|^{p_l+1-\iota}+|y|^{p_h+1-\iota})}{|x|^{1-\iota}}1_{\{|x|\neq 0\}}
\end{eqnarray*}
for some $0<\iota<1$.

Substituting in the above inequality $x$ with
$\Delta_n^{-1/\beta}\int_{(i-1)\Delta_n}^{i\Delta_n}\overline{\sigma}_{s-}\,dL_s$
and $y$ with $a_i\Delta_n^{-1/\beta}$ and using the fact that
$(p|x|^{p-1}-q|x|^{q-1})\operatorname{sign}\{x\}1_{\{|x|\neq 0, 2|y|\leq
|x|\}}$ is odd in $x$, we get
\begin{eqnarray*}\label{eq:at5}
&&
\mathbb{E}\bigl(\mathbb{E}_{i-1}^n\bigl([c_i(q)-c_i(p)]1_{\{C_i^n\}}\bigr)\bigr)^2
\\
&&\qquad\leq
K(p-q)^2\bigl(|\Delta_n|^{2(p_l+1-\iota)(1-1/\beta)}+|\Delta_n|^{2(p_h+1-\iota)(1-1/\beta)}\bigr)
\end{eqnarray*}
for some $\iota<p_l-\frac{2-\beta}{2(\beta-1)}$. For
$\widetilde{A}_{32}^n(p,q)$ we have for sufficiently small $\Delta_n$
\[
\sup_{p,q\in[p_l,p_h]} \widetilde{A}_{32}^n(p,q)\leq K\Delta_n^{1/2}\sum_{i=1}^{[T/\Delta_n]}a_i^{p_l}\Delta_n^{-p_l/\beta}1_{\{(C_i^{n})^c\}}.
\]
Then using the definition of the set $(C_i^{n})^c$ and the calculation
in \refeq{eq:an_3a} we can conclude
%e6.22 ###
\begin{equation}
\limsup_{\Delta_n\downarrow 0}\mathbb{E}\Bigl(\sup_{p,q\in[p_l,p_h]}\widetilde{A}_{32}^n(p,q)\Bigr)=0.
\end{equation}
Combining the above results we get the tightness of
$\widehat{V}_T(q,X,\Delta_n)$ on the space of continuous functions of
$p$ in the interval $[p_l,p_h]$.
\end{pf}

%s6.2 ###
\subsection{\texorpdfstring{Proof of Remark \protect\ref{rema38}}{Proof of Remark 3.8.}}

In what follows we denote
\[
\chi_i^n:=\Delta_n^{p/\beta}\bigl(|\Delta_n^{-1/\beta}\Delta_i^n X|^p-\Pi_{A,\beta}^{\beta}\mu_p(\beta)\bigr).
\]
It is no restriction, of
course, to assume that the constant $A$ in \refeq{eq:nu_1_2}
corresponds to that of a standard stable, and we proceed in the proofs
with that assumption. In view of Theorem XVII.2.2 in \cite{Feller71} we
need to prove the following:
%e6.24 ###
%e6.23 ###
\begin{eqnarray}
\label{eq:st_1}\frac{1}{\Delta_n}\mathbb{E}\bigl(\chi_i^n 1_{\{|\chi_i^n|\leq 1\}}\bigr)
&\rightarrow&
-2\frac{\beta}{\beta-p}\frac{A}{\beta},
\\
\label{eq:st_2}\qquad\frac{1}{\Delta_n}\bigl[\mathbb{E}\bigl((\chi_i^n)^21_{\{|\chi_i^n|\leq K\}}\bigr)-\bigl(\mathbb{E}\bigl(\chi_i^n1_{\{|\chi_i^n|\leq K\}}\bigr)\bigr)^2\bigr]
&\rightarrow&
2K^{2-\beta/p}\frac{\beta}{2p-\beta}\frac{A}{\beta},\hspace*{-8pt}
\\
\label{eq:st_3}\frac{1}{\Delta_n}\mathbb{E}\bigl(1_{\{\chi_i^n>K\}}\bigr)
&\rightarrow&
2K^{\beta/p}\frac{A}{\beta}\quad\mbox{and}\nonumber
\\[-8pt]\\[-8pt]
\frac{1}{\Delta_n}\mathbb{E}\bigl(1_{\{\chi_i^n<-K\}}\bigr)
&\rightarrow&
0,\nonumber
\end{eqnarray}
where $K>0$ is an arbitrary positive constant.

We recall that $X$ is symmetric stable process plus a drift, that is,
$X_t\stackrel{d}{=}L_t+a t$, where $L_t$ denotes symmetric stable
process with \Lvb density equal to $\nu_1(x)$ in \refeq{eq:nu_1_2} and
$a=m_d+\int_{\mathbb{R}}(x-\kappa(x))\nu_1(x)\,dx$ when \mbox{$\beta>1$} and
$a=0$ when \mbox{$\beta\leq 1$}. Using the self-similarity of the symmetric
stable we have $\Delta_n^{-1/\beta}\Delta_i^n
X\stackrel{d}{=}L_1+a\Delta_n^{1-1/\beta}$.

First we state several basic facts about the stable distribution that
we make use of in the proof. We recall that for the tail of the
symmetric stable we have (see, e.g., \cite{Zolotarev86})
$\mathbb{P}(L_1>x)\sim \mathbb{P}(L_1<-x)\sim
\frac{A}{\beta}\frac{1}{x^{\beta}}$ as $x\uparrow +\infty$ where for
two\vspace*{-2pt} functions $f(\cdot)$ and $g(\cdot)$, $f(\Delta_n)\sim g(\Delta_n)$
means $\lim_{\Delta_n\downarrow
0}\frac{f(\Delta_n)}{g(\Delta_n)}$$=$$1$. Therefore the tail
probability of the stable distribution varies regularly at infinity,
and we can use this fact and Theorems 8.1.2 and 8.1.4 in
\cite{Teugels87} to write for $p\in(\beta/2,\beta)$
%e6.25 ###
\begin{eqnarray}
\label{eq:st_4}\mathbb{E}\bigl(|L_1|^p1_{\{L_1>x\}}\bigr)
&\sim&
\mathbb{E}\bigl(|L_1|^p1_{\{L_1<-x\}}\bigr)\nonumber
\\[-8pt]\\[-8pt]
&\sim&
x^{p-\beta}\frac{\beta}{\beta-p}\frac{A}{\beta},\nonumber
\\
\label{eq:st_5}\mathbb{E}\bigl(|L_1|^{2p}1_{\{|L_1|\leq x\}}\bigr)
&\sim&
2x^{2p-\beta}\frac{\beta}{2p-\beta}\frac{A}{\beta},
\end{eqnarray}
as $x\uparrow\infty$. We continue with the proof of
\refeq{eq:st_1}--\refeq{eq:st_3}. We start with showing \refeq{eq:st_1}.
First we have
%e6.26 ###
\begin{eqnarray}\label{eq:st_6}
\frac{1}{\Delta_n}\mathbb{E}(\chi_i^n)
&=&
\Delta_n^{p/\beta-1}\mathbb{E}(|L_1+a\Delta_n^{1-1/\beta}|^p-|L_1|^p)\nonumber
\\
&&{}+
\Delta_n^{p/\beta-1}\mathbb{E}\bigl(|L_1|^p-\Pi_{A,\beta}^{\beta}\mu_p(\beta)\bigr)
\\
&\rightarrow&
0.\nonumber
\end{eqnarray}
We note that the second term on the right-hand side of \refeq{eq:st_6}
is identically zero, while the convergence of the first term can be
split into two cases. First, when $p\leq 1$ the result follows from the
bound for the term $A_{i2}^n$ in \refeq{eq:an_3} and \refeq{eq:an_3a} in
the proof of Theorem~\ref{thms:asl}, provided $p>1/\beta$. When $p>1$
the convergence follows from a trivial application of the Taylor
expansion.

Second using the rate of decay of the tail probability of the stable
distribution we have
\[%\label{eq:st_7}
\Delta_n^{p/\beta-1}\mathbb{P}\bigl(\bigl||L_1+a\Delta_n^{1-1/\beta}|^p-\Pi_{A,\beta}^{\beta}\mu_p(\beta)\bigr|>\Delta_n^{-p/\beta}\bigr)
\rightarrow 0.
\]
Third using a Taylor expansion around $L_1$ and the fact that we
evaluate $L_1$ on a set growing to infinity at the rate
$\Delta_n^{-1/\beta}$,  we have
\[%\label{eq:st_8}
\Delta_n^{p/\beta-1}\mathbb{E}(|L_1+a\Delta_n^{1-1/\beta}|^p-|L_1|^p)
1_{\{||L_1+a\Delta_n^{1-1/\beta}|^p-\Pi_{A,\beta}^{\beta}\mu_p(\beta)|>\Delta_n^{-p/\beta}\}}
\rightarrow 0.
\]
Thus to prove \refeq{eq:st_1} we need to show
\[%\label{eq:st_9}
\Delta_n^{p/\beta-1}\mathbb{E}|L_1|^p1_{\{||L_1+a\Delta_n^{1-1/\beta}|^p-\Pi_{A,\beta}^{\beta}\mu_p(\beta)|>\Delta_n^{-p/\beta}\}}
\rightarrow 2\frac{\beta}{\beta-p}\frac{A}{\beta}.
\]
But this follows from \refeq{eq:st_4} with
\[
x=\bigl(\bigl(\Pi_{A,\beta}^{\beta}\mu_p(\beta)+\Delta_n^{-p/\beta}\bigr)^{1/p}\pm a\Delta_n^{1-1/\beta}\bigr),
\]
and hence we are done. We turn now to
\refeq{eq:st_2}. It is easy to show that
\begin{eqnarray*}%\label{eq:st_10}
&&
\Delta_n^{2p/\beta-1}\mathbb{E}(|L_1+a\Delta_n^{1-1/\beta}|^{2p}-|L_1|^{2p})
1_{\{||L_1+a\Delta_n^{1-1/\beta}|^p-\Pi_{A,\beta}^{\beta}\mu_p(\beta)|\leq
K\Delta_n^{-p/\beta}\}}
\\
&&\qquad\rightarrow
0.
\end{eqnarray*}
Therefore, \refeq{eq:st_2} will follow if we can show
\begin{eqnarray}\label{eq:st_11}
&&
\Delta_n^{2p/\beta-1}\mathbb{E}|L_1|^{2p}1_{\{||L_1+a\Delta_n^{1-1/\beta}|^p-\Pi_{A,\beta}^{\beta}\mu_p(\beta)|\leq K\Delta_n^{-p/\beta}\}}\nonumber
\\[-8pt]\\[-8pt]
&&\qquad\rightarrow
2K^{2-\beta/p}\frac{\beta}{2p-\beta}\frac{A}{\beta}.\nonumber
\end{eqnarray}
To show \refeq{eq:st_11} we can apply \refeq{eq:st_5} with
\[
x=\bigl(\bigl(\Pi_{A,\beta}^{\beta}\mu_p(\beta)+K\Delta_n^{-p/\beta}\bigr)^{1/p}\pm
a\Delta_n^{1-1/\beta}\bigr).
\]

Finally, \refeq{eq:st_3} follows trivially from the expression for the
tail probability of a stable stated earlier.

\subsection{\texorpdfstring{Proof of Corollary~\protect\ref{thms:asp}}{Proof of Corollary 4.1.}}

Again, as in the proof of Theorem~\ref{thms:asl} we will show only part
(b), the proof of part (a) being identical. Since the process $X$ has
no fixed time of discontinuity, the result of Lemma~\ref{lema_fd}
implies that the convergence in \refeq{eq:js_0} holds for an arbitrary
\textit{fixed} $T>0$. Then, there is a set $\Omega_n$ on which
$2V_T(X,p,2\Delta_n)\neq V_T(X,p,\Delta_n)$ for $p\in[p_l,p_h]$ and
from Theorem~\ref{thms:cp_gen} (under the conditions of this theorem)
$\Omega_n \rightarrow\Omega$. On $\Omega_n$ $b_{X,T}(p)$ is a
continuous transformation of $V_T(X,p,2\Delta_n)$ and
$V_T(X,p,\Delta_n)$, and thus Lemma~\ref{lema_fd} implies the
finite-dimensional convergence of the sequences on the left-hand sides
of \refeq{eq:asp_a} and \refeq{eq:asp_b}. Similarly, since tightness is
preserved under continuous transformations, using Lemma~\ref{lema_ti}
we have that the left-hand sides of \refeq{eq:asp_a} and
\refeq{eq:asp_b} are tight. Hence the result of Theorem~\ref{thms:asp}
follows.

%s6.4 ###
\subsection{\texorpdfstring{Proof of Theorem~\protect\ref{thms:ase}}{Proof of Theorem 4.1.}}

We first show the result for the case when $w(u)$ is continuous on
$[\tau_1^{*},\tau_2^{*}]$. Set
\[
\tau_1(z)=f_l(z)\quad\mbox{and}\quad\tau_2(z)=f_h(z).
\]
Since $\tau_1(z)$ is continuous in a neighborhood of $\beta_{X,T}$ and
$\tau_1(\beta_{X,T})>\frac{\beta'}{2-\beta'}$ as well as
$\tau_2(\beta_{X,T})<\beta_{X,T}/2$ when $X$ is given by
\refeq{eq:X_c}, then there are $z_{*}<\beta_{X,T}<z^{*}$ such that for
all $z\in(z_{*},z^{*})$ $\Rightarrow$
$\tau_1(z)>\frac{\beta'}{2-\beta'}$ and $\tau_2(z)<\beta_{X,T}/2$.
Similarly if $X$ is given by \refeq{eq:X_d}, then
$\beta_{X,T}\equiv\beta$, and due to the assumptions of the theorem,
there exist $z_{*}<\beta<z^{*}$ such that for $z\in(z_{*},z^{*})$
$\Rightarrow$ $\tau_1(z)>(\frac{2-\beta}{2(\beta-1)}\vee
\frac{\beta\beta'}{2(\beta-\beta')})$ and $\tau_2(z)<\beta/2$ when
$\beta>\sqrt{2}$ and $z\in(z_{*},z^{*})$ $\Rightarrow$
$\tau_1(z)>\frac{\beta\beta'}{2(\beta-\beta')}$ and $\tau_2(z)<\beta/2$
when $\beta\leq \sqrt{2}$.

Denote with $A$ the subset of $(z_{*},z^{*})$ for which $\tau_1(z)$ and
$\tau_2(z)$ are continuously differentiable. From the assumptions of
Theorem~\ref{thms:ase} the set $A$ contains a neighborhood of
$\beta_{X,T}$. Then, using\vadjust{\goodbreak} a Taylor expansion on the set
$B_n:=\{\omega\dvtx \hat{\beta}_{X,T}^{fs}\in A\cap\Omega_n\}$ where
$\Omega_n$ is the set defined in the proof of Corollary~\ref{thms:asp}
above, we can write
%e6.27 ###
\begin{eqnarray}\label{eq:ase1}
\Delta_n^{-1/2}(\hat{\beta}_{X,T}^{ts}-\beta_{X,T})
&=&
1_{B_n}\int_{\tau_1^{*}}^{\tau_2^{*}}w(u)\bigl\{\Delta_n^{-1/2}\bigl(b_{X,t}(u)-\beta_{X,T}\bigr)\bigr\}\,du\nonumber
\\
&&{}+
1_{B_n}\Delta_n^{-1/2}\Theta_T(\overline{\beta}_{X,T})(\hat{\beta}_{X,T}^{fs}-\beta_{X,T})
\\
&&{}+
1_{B_n^c}\Delta_n^{-1/2}(\hat{\beta}_{X,T}^{ts}-\beta_{X,T}),\nonumber
\end{eqnarray}
where $\overline{\beta}_{X,T}$ is between $\hat{\beta}_{X,T}^{fs}$
and $\beta_{X,T}$ and
\begin{eqnarray*}%\label{eq:ase2}
\Theta_T(z)
&=&
w(\tau_2(z))\nabla_z\tau_2(z)\bigl(b_{X,T}(\tau_2(z))-\beta_{X,T}\bigr)
\\
&&{}-
w(\tau_1(z))\nabla_z\tau_1(z)\bigl(b_{X,T}(\tau_1(z))-\beta_{X,T}\bigr).
\end{eqnarray*}
The last term on the right-hand side of \refeq{eq:ase1} is
asymptotically negligible because $\hat{\beta}_{X,T}^{fs}$ is
consistent for $\beta_{X,T}$. We now show that the\vspace*{-3pt} second term in
\refeq{eq:ase1} is asymptotically negligible. First note that\vspace*{-2pt} since
$\hat{\beta}_{X,T}^{fs}\stackrel{\mathbb{P}}{\longrightarrow}\beta_{X,T}$
we also have
$\overline{\beta}_{X,T}\stackrel{\mathbb{P}}{\longrightarrow}\beta_{X,T}$.
Then to establish the asymptotic negligibility it  suffices to show
that
%e6.28 ###
\begin{equation}\label{eq:ase3}
\qquad\mathbb{P}\biggl(\Delta_n^{-1/2}\int_{\overline{\tau}_1}^{\overline{\tau}_2}\bigl|\bigl(b_{X,T}(u)-\beta_{X,T}\bigr)w(u)\bigl|\,du>\varepsilon\biggr)\downarrow
0\qquad\mbox{for }\varepsilon\uparrow+\infty,
\end{equation}
where $\overline{\tau}_1:=\tau_1(\overline{\beta}_{X,T})$ and
$\overline{\tau}_2:=\tau_2(\overline{\beta}_{X,T})$. For any
$\varepsilon>0$ we have
\begin{eqnarray*}%\label{eq:ase4}
&&
\mathbb{P}\biggl(\Delta_n^{-1/2}\int_{\overline{\tau}_1}^{\overline{\tau}_2}\bigl|\bigl(b_{X,T}(u)-\beta_{X,T}\bigr)w(u)\bigr|\,du>\varepsilon\biggr)
\\
&&\qquad\leq
\mathbb{P}(\overline{\beta}_{X,T}\in A^c)
+
\mathbb{P}\biggl(1_{\{\overline{\beta}_{X,T}\in A\}}\int_{\overline{\tau}_1}^{\overline{\tau}_2}\bigl|\Delta_n^{-1/2}\bigl(b_{X,T}(u)-\beta_{X,T}\bigr)w(u)\bigr|\,du>\varepsilon\biggr).
\end{eqnarray*}
The first probability in the second line of \refeq{eq:ase3} is
converging to $0$ as $\Delta_n\downarrow 0$, while the second one
converges to zero as $\varepsilon\uparrow +\infty$. This is because when
$\overline{\beta}_{X,T}\in A$, $\overline{\tau}_1>p_l$ and
$\overline{\tau}_2<p_h$ where $p_l<p_h$ are some constants that satisfy
the conditions of  Theorem~\ref{thms:asl}, and as a consequence of this
theorem $\Delta_n^{-1/2}(b_{X,t}(u)-\beta_{X,T})$ converges uniformly
in $u$ for $u\in[p_l,p_h]$.

We are left with the first term in \refeq{eq:ase1}. Using the uniform
convergence result of Theorem~\ref{thms:asl}, the fact that the
integration over a bounded interval is continuous for the uniform
metric on the space of continuous functions (in fact for this even
finite  dimensional convergence suffices) we have
\[
\int_{\tau_1^{*}}^{\tau_2^{*}}w(u)\bigl\{\Delta_n^{-1/2}\bigl(b_{X,T}(u)-\beta_{X,T}\bigr)\bigr\}\,du
\stackrel{\mathcal{L}-s}{\longrightarrow}
\int_{\tau_1^{*}}^{\tau_2^{*}}Z_{\beta_{X,T}}(u)w(u)\,du,
\]
where
\[
Z_{\beta}(u)=
\cases{
\displaystyle\frac{\beta^2}{u\ln 2}\frac{1}{T\Pi_{A,\beta}^{u/\beta}\mu_{u}(\beta)}\bigl(\Psi_{\beta,T}^{(2)}(u)-2^{1-u/\beta}\Psi_{\beta,T}^{(1)}(u)\bigr)
&\quad if $\beta<2$,\cr
\displaystyle\frac{4}{u\ln 2}\frac{1}{T|\sigma|^u\mu_{u}(2)}\bigl(\Psi_{2,T}^{(2)}(u)-2^{1-u/2}\Psi_{2,T}^{(1)}(u)\bigr)
&\quad if $\beta=2$,
}
\]
and $\Psi_{\beta,T}^{(1)}$ and $\Psi_{\beta,T}^{(2)}$ are the first and
second elements, respectively, of the limiting Gaussian process of part
(a) and (b) of Theorem~\ref{thms:asl}. The proof of
Theorem~\ref{thms:ase} for the case of continuous $w(u)$ then easily
follows. The proof in the case of $w(u)$ being Dirac mass at some point
follows from the proof of Corollary~\ref{thms:tse} given below.

%s6.5 ###
\subsection{\texorpdfstring{Proof of Corollary~\protect\ref{thms:tse}}{Proof of Corollary 4.2.}}

Denote with $A$ the set of values of $z$ for which $f(z)\in(p_l,p_h)$
for some $0<p_l<p_h<\beta_{X,T}/2$ satisfying the conditions of
Theorem~\ref{thms:asl} in the different cases for $\beta_{X,T}$.
Finally, set $B_n:=\{\omega\dvtx\hat{\beta}_{X,T}^{fs}\in A\cap
\Omega_n\}$. We know that this set contains neighborhood of
$\beta_{X,T}$ because of the continuity of $f(\cdot)$ and the fact that
$p^*\in (p_l,p_h)$. Then we can write
\begin{eqnarray*}
\Delta_n^{-1/2}\bigl(\hat{\beta}_{X,T}^{ts}-\beta_{X,T}\bigr)
&=&
1_{B_n}\Delta_n^{-1/2}\bigl(b_{X,T}(\tau^{*})-\beta_{X,T}\bigr)+1_{B_n^c}\Delta_n^{-1/2}(\hat{\beta}_{X,T}^{ts}-\beta_{X,T})
\\
&&{}+
1_{B_n}\Delta_n^{-1/2}\Theta_T(\overline{f(\hat{\beta}_{X,T}^{fs})})\bigl(f(\hat{\beta}_{X,T}^{fs})-f(\beta_{X,T})\bigr),
\end{eqnarray*}
where $\overline{f(\hat{\beta}_{X,T}^{fs})}$ is between
$f(\hat{\beta}_{X,T}^{fs})$ and $f(\beta_{X,T})$ and
\begin{eqnarray*}
\Theta_T(z)
&=&
\Theta_T^{(1)}(z)+\Theta_T^{(2)}(z),
\qquad
\Theta_T^{(1)}(z)
=
\frac{b_{X,T}(z)-\beta_{X,T}}{z}-\frac{b_{X,T}^2(z)-\beta^2_{X,T}}{\beta_{X,T}z},
\\
\Theta_T^{(2)}(z)
&=&
\frac{b_{X,T}^2(z)}{z\ln 2}\biggl(\frac{\nabla_z
[\Delta_n^{1-z/\beta_{X,T}}V_T(z,X,\Delta_n)]}{\Delta_n^{1-z/\beta_{X,T}}V_T(z,X,\Delta_n)}
\\
&&\hphantom{\frac{b_{X,T}^2(z)}{z\ln 2}\biggl(}{}-
\frac{\nabla_z[(2\Delta_n)^{1-z/\beta_{X,T}}V_T(z,X,2\Delta_n)]}{(2\Delta_n)^{1-z/\beta_{X,T}}V_T(z,X,2\Delta_n)}\biggr).
\end{eqnarray*}
The result of Corollary~\ref{thms:tse} then will follow if we can show
that
$\Delta_n^{-1/2}\Theta_T(\overline{f(\hat{\beta}_{X,t}^{fs})})$ is
bounded in probability on the set $B_n$. But this holds true because we
can prove exactly as in Theorem~\ref{thms:asl} that
\begin{eqnarray*}
&&
\Delta_n^{-1/2}\Biggl(\Delta_n\sum_{i=1}^{[T/\Delta_n]}|\Delta_n^{-1/\beta_{X,T}}\Delta_i^n X|^p\ln|\Delta_n^{-1/\beta_{X,T}}\Delta_i^n X|1_{\{|\Delta_i^n X|>0\}}
\\
&&\hspace*{189pt}{}-
T\mathbb{E}(|L_1|^p\ln|L_1|)\Biggr)
\end{eqnarray*}
converges uniformly in $p$ (under the same conditions for the power as
in that theorem).

%s6.6 ###
\subsection{\texorpdfstring{Proof of Theorem~\protect\ref{thms:sv}}{Proof of Theorem 3.4.}}
We do not show here part (a). The finite-dimensional convergence for
this case (without jumps in $X$) has been already shown in
\cite{BNGJPS05} (extending their result to the case with jumps
satisfying the conditions of Theorem~\ref{thms:sv}, part(a) follows
trivially using the subadditivity of $|x|^p$ for $p\leq 1$). The
tightness can be shown in exactly the same way as part (b) (i.e., in
the decomposition in equation (8.2); in \cite{BNGJPS05} we can apply
the same techniques as in the proof of our Lemma~\ref{lema_ti}).

\begin{pf*}{Proof of part(b)}  We will establish only the
finite-dimensional convergence, the proof the tightness is done exactly
as in Lemma~\ref{lema_ti}. Also we will prove the finite-dimensional
convergence for a fixed $p$ and the second element of the vector on the
left-hand side in \refeq{eq:sv_b}. The generalization will follow
immediately.

As in the previous proofs we assume that $A$ in \refeq{eq:nu_1_2}
corresponds to that of a standard stable. Upon using a localization
argument as in \cite{J06b} we can and will assume the following
stronger assumption on the various processes in
\refeq{eq:X_d_sv} and \refeq{eq:X_d_sv_sv}:

\textit{We have
$|m_{dt}|+|b_t|+|\sigma_{2t}|+|\sigma_{2t}|^{-1}+|\widetilde{\sigma}_{2t}|\leq
K$ and $|\delta(t,\mathbf{x})|\leq \gamma(\mathbf{x})\leq K$ for some
positive constant $K$ which bounds also the coefficients in the \Ito
semimartingale representations of the processes $m_{dt}$ and
$\widetilde{\sigma}_{2t}$; $\int_{\mathbb{R}}1_{|x|>K}\nu(x)\,\allowbreak dx=0$.}

We can make the following decomposition:
\[
\Delta_n^{-1/2}\biggl(\Delta_n^{1-p/\beta}V_T(X,p,\Delta_n)-\mu_p(\beta)\int_0^T|\sigma_{2s}|^p\,ds\biggr)=\sum_{i=1}^5A_i,
\]
%e6.30 ###
%e6.29 ###
\begin{eqnarray}
A_1
&=&
\Delta_n^{1/2}\sum_{i=1}^{[T/\Delta_n]}\biggl(|\Delta_n^{-1/\beta}\Delta_i^n
\overline{X}|^p-\mu_p(\beta)\biggl(\frac{1}{\Delta_n}\int_{(i-1)\Delta_n}^{i\Delta_n}|\overline{\sigma}_{2s}|^{\beta}\,ds\biggr)^{p/\beta}\biggr),\nonumber
\\
A_2
&=&
\mu_p(\beta)\Delta_n^{1/2}\sum_{i=1}^{[T/\Delta_n]}a_{i2},\nonumber
\\
\eqntext{
a_{i2}=\biggl(\frac{1}{\Delta_n}\int_{(i-1)\Delta_n}^{i\Delta_n}|\overline{\sigma}_{2s}|^{\beta}\,ds\biggr)^{p/\beta}
-
\bigl|\sigma_{2,(i-1)\Delta_n}\bigr|^p,}
\\
A_3
&=&
\mu_p(\beta)\Delta_n^{1/2}\sum_{i=1}^{[T/\Delta_n]}a_{i3},\nonumber
\\
\eqntext{
a_{i3}=\frac{1}{\Delta_n}\int_{(i-1)\Delta_n}^{i\Delta_n}\bigl(\bigl|\sigma_{2,(i-1)\Delta_n}\bigr|^p
-
|\overline{\sigma}_{2s}|^p\bigr)\,ds,}
\\
A_4
&=&
\Delta_n^{-1/2}\mu_p(\beta)\int_0^T(|\overline{\sigma}_{2s}|^p-|\sigma_{2s}|^p)\,ds,\nonumber
\\
A_5
&=&
\Delta_n^{1/2}\sum_{i=1}^{[T/\Delta_n]}(|\Delta_n^{-1/\beta}\Delta_i^nX|^p-|\Delta_n^{-1/\beta}\Delta_i^n\overline{X}|^p),\nonumber
\end{eqnarray}
where for $i=1,\ldots,[T/\Delta_n]$ and $s\in [(i-1)\Delta_n,i\Delta_n)$
\begin{eqnarray*}
\overline{\sigma}_{2s}
&=&
\sigma_{2,(i-1)\Delta_n}+\widetilde{\sigma}_{2,(i-1)\Delta_n}\bigl(W_s-W_{(i-1)\Delta_n}\bigr),
\\
\overline{X}_{s}
&=&
X_{(i-1)\Delta_n}+\int_{(i-1)\Delta_n}^sm_{d,(i-1)\Delta_n}\,du+\int_{(i-1)\Delta_n}^s\int_{\mathbb{R}}\overline{\sigma}_{2u-}\kappa(x)\widetilde{\mu}(du,dx)
\\
&&{}+
\int_{(i-1)\Delta_n}^s\int_{\mathbb{R}}\overline{\sigma}_{2u-}\kappa'(x)\mu(du,dx),\qquad s\in \bigl[(i-1)\Delta_n,i\Delta_n\bigr).
\end{eqnarray*}

We start with $A_1$. We can apply directly Lemma~\ref{lema_fd} to show
that $A_1$ converges stably to the limit on the right-hand side of
\refeq{eq:sv_b} (recall our stronger assumption on the process
$\sigma_2$ stated at the beginning of the proof). We continue with the
term $A_2$ which we now show is asymptotically negligible. First we
denote the set
\[
B_{i,n}:=\Bigl\{\omega\dvtx\sup_{s\in[(i-1)\Delta_n,i\Delta_n]}\bigl|\sigma_{2,(i-1)\Delta_n}-\overline{\sigma}_{2s}\bigr|
>
0.5\sigma_{2,(i-1)\Delta_n}\Bigr\}.
\]
Then, using the exponential inequality for continuous martingales with
bounded variation [see, e.g., \cite{RY}] it is easy to derive
\[
\bigl|\mathbb{E}_{i-1}^n1_{\{B_{i,n}\}}a_{i2}^n\bigr|
\leq
K e^{-K/\Delta_n},
\qquad
\mathbb{E}_{i-1}^n1_{\{B_{i,n}\}}(a_{i2}^n)^2
\leq
K e^{-K/\Delta_n}.
\]
Using a second-order Taylor expansion and the fact that
$\overline{\sigma}_{2s}$ is bounded from below on the set
$(B_{i,n})^c$, we get
\begin{eqnarray*}
&&
\bigl|\mathbb{E}_{i-1}^n\bigl(1_{\{(B_{i,n})^c\}}a_{i2}^n\bigr)\bigr|
\\
&&\qquad\leq
K\mathbb{E}_{i-1}^n\biggl(1_{\{(B_{i,n})^c\}}\frac{1}{\Delta_n}\int_{(i-1)\Delta_n}^{i\Delta_n}\bigl||\overline{\sigma}_{2s}|^{\beta}-\bigl|\sigma_{2,(i-1)\Delta_n}\bigr|^{\beta}\bigr|\,ds\biggr)^2
\\
&&\qquad\quad{}+
K\mathbb{E}_{i-1}^n\biggl(\frac{1}{\Delta_n}\int_{(i-1)\Delta_n}^{i\Delta_n}\bigl(\overline{\sigma}_{2s}-\overline{\sigma}_{2,(i-1)\Delta_n}\bigr)^2\,ds\biggr)
\\
&&\qquad\leq
K\Delta_n,
\end{eqnarray*}
where we also made use of the following inequality:
\begin{eqnarray*}
&&
\biggl|\mathbb{E}_{i-1}^n\biggl(1_{\{(B_{i,n})^c\}}\int_{(i-1)\Delta_n}^{i\Delta_n}\bigl(\overline{\sigma}_{2s}-\sigma_{2,(i-1)\Delta_n}\bigr)\,ds\biggr)\biggr|
\\
&&\qquad=
\biggl|\mathbb{E}_{i-1}^n\biggl(1_{\{B_{i,n}\}}\int_{(i-1)\Delta_n}^{i\Delta_n}\bigl(\overline{\sigma}_{2s}-\sigma_{2,(i-1)\Delta_n}\bigr)\,ds\biggr)\biggr|\leq
K e^{-K/\Delta_n}.
\end{eqnarray*}
Finally, a first-order Taylor expansion together with the fact that
$\overline{\sigma}_{2s}$ is bounded from below on the set $(B_{i,n})^c$
gives
\begin{eqnarray*}
&&
\mathbb{E}_{i-1}^n\bigl(1_{\{(B_{i,n})^c\}}a_{i2}^n\bigr)^2
\\
&&\qquad\leq
K\mathbb{E}_{i-1}^n\biggl(1_{\{(B_{i,n})^c\}}\frac{1}{\Delta_n}\int_{(i-1)\Delta_n}^{i\Delta_n}\bigl||\overline{\sigma}_{2s}|^{\beta}-\bigl|\sigma_{2,(i-1)\Delta_n}\bigr|^{\beta}\bigr|\,ds\biggr)^2
\\
&&\qquad\leq
K\Delta_n.
\end{eqnarray*}
Combining the above two inequalities we get
\[
\cases{
\displaystyle\Delta_n^{1/2}\sum_{i=1}^{[T/\Delta_n]}\mathbb{E}_{i-1}^na_{i2}\stackrel{\mathrm{u.c.p.}}{\longrightarrow}0,\cr
\displaystyle\Delta_n\sum_{i=1}^{[T/\Delta_n]}\mathbb{E}_{i-1}^n(a_{i2})^2\stackrel{\mathrm{u.c.p.}}{\longrightarrow}0.
}
\]
This implies the asymptotic negligibility of $A_2$. We continue with
$A_{3}$.  We can use the standard inequality $|a+b|^p\leq |a|^p+|b|^p$
for $0<p\leq 1$ as well as H\"{o}lder's inequality to get
\[
\bigl|\mathbb{E}_{i-1}^n1_{\{B_{i,n}\}}a_{i3}^n\bigr|\leq K e^{-K/\Delta_n},
\qquad\mathbb{E}_{i-1}^n1_{\{B_{i,n}\}}(a_{i3}^n)^2
\leq
K e^{-K/\Delta_n}.
\]
Similar inequalities as for $a_{i2}$ on the set $(B_{i,n})^c$ give
\[
\bigl|\mathbb{E}_{i-1}^n1_{\{(B_{i,n})^c\}}a_{i3}^n\bigr|\leq K\Delta_n,
\qquad\mathbb{E}_{i-1}^n\bigl(1_{\{(B_{i,n})^c\}}a_{i3}^n\bigr)^2
\leq
K\Delta_n.
\]
These two inequalities establish the asymptotic negligibility of $A_3$.
We continue with $A_4$. First, for some $\varepsilon>0$ denote the set
$B_{i,\varepsilon}^n:=\{\omega\dvtx\sup_{s\in[(i-1)\Delta_n,i\Delta_n]}|\allowbreak\sigma_{2s}-\overline{\sigma}_{2s}|>\varepsilon\}$.
Then we can decompose $A_4$ into
\begin{eqnarray*}%\label{eq:sv_a2_1}
A_4
&=&
\mu_p(\beta)(C_1+C_2+C_3),
\\
C_1
&=&
\Delta_n^{-1/2}\sum_{i=1}^{[T/\Delta_n]}\int_{(i-1)\Delta_n}^{i\Delta_n}g(\sigma_{2s})(\sigma_{2s}-\overline{\sigma}_{2s})\,ds,
\\
C_2
&=&
p(p-1)\Delta_n^{-1/2}\sum_{i=1}^{[T/\Delta_n]}1_{(B_{i,\varepsilon}^n)^c}\int_{(i-1)\Delta_n}^{i\Delta_n}|\sigma_{2s}^{*}|^{p-2}(\sigma_{2s}-\overline{\sigma}_{2s})^2\,ds,
\\
C_3
&=&
\Delta_n^{-1/2}\sum_{i=1}^{[T/\Delta_n]}1_{B_{i,\varepsilon}^n}\int_{(i-1)\Delta_n}^{i\Delta_n}\bigl(|\overline{\sigma}_{2s}|^p-|\sigma_{2s}|^p-(\sigma_{2s}-\overline{\sigma}_{2s})g(\sigma_{2s})\bigr)\,ds,
\end{eqnarray*}
where $\sigma_{2s}^{*}$ is a number between $\sigma_{2s}$ and
$\overline{\sigma}_{2s}$ and $g(x)=p\mbox{
}\operatorname{sign}\{x\}|x|^{p-1}$. Note that for $\varepsilon$ sufficiently
small $C_2$ is well defined because of the boundedness from below of
$|\sigma_{2s}|$.

Using the Burkholder--Davis--Gundy inequality, H\"{o}lder's inequality,
the assumption of \Ito semimartingale for the process
$\widetilde{\sigma}_2$ (due to which the leading term in
$\sigma_{2s}-\overline{\sigma}_{2s}$ is
$\int_{(i-1)\Delta_n}^s\int_{\mathbb{R}^2}\kappa(\delta(u,\mathbf{x}))\tilde{\underline{\mu}}(u,\mathbf{x})$)
and the integrability condition for the dominating function of the
jumps in $\sigma_{2t}$, $\gamma(\mathbf{x})$, in \refeq{eq:X_d_sv_a},
we have for $s\in[(i-1)\Delta_n,i\Delta_n)$
%e6.31 ###
\begin{equation}\label{eq:sv_a2_6}
\cases{
\mathbb{E}_{i-1}^n|\sigma_{2s}-\overline{\sigma}_{2s}|^{p}\leq K|s-(i-1)\Delta_n|^{p/\beta-\varepsilon}\cr
\qquad\mbox{for $p\leq \beta$, $\forall \varepsilon>0$},\cr
\mathbb{E}_{i-1}^n|\sigma_{2s}-\overline{\sigma}_{2s}|^{p}\leq K|s-(i-1)\Delta_n|\cr
\qquad\mbox{for $p> \beta$}
}
\end{equation}
for some constant $K$ that does not depend on $\Delta_n$. We will show
that the three terms $C_1$, $C_2$ and $C_3$ are asymptotically
negligible. For $C_1$ and $C_2$ we make use of the fact that a
sufficient condition for asymptotic negligibility\vspace*{-2pt} of
$\sum_{i=1}^{[T/\Delta_n]}\xi_i^n$, where $\xi_i^n$ is
$\mathscr{F}_{i\Delta_n}$-measurable, is
$\sum_{i=1}^{[T/\Delta_n]}\mathbb{E}_{i-1}^n|\xi_i^n|\stackrel{\mathbb{P}}{\longrightarrow}0$
(see Theorem VIII.2.27 of \cite{JS} (or the first part of Lemma 4.1 in
\cite{J07})). Note that for $C_2$ we use the fact that
$\sigma_{2s}^{*}$ is bounded by a constant on the set
$(B_{i,\varepsilon}^n)^c$. For $C_3$ we can first make use of Doob's
inequality to show that $\mathbb{P}(\omega\in B_{i,\varepsilon}^n)\leq
K\Delta_n$ for\vspace*{-2pt} some constant $K$ that depends on~$\varepsilon$. Then,
since $\mathbb{E}(\int_{(i-1)\Delta_n}^{i\Delta_n}
(|\sigma_{2s}|^p-|\overline{\sigma}_{2s}|^p-(\sigma_{2s}-
\overline{\sigma}_{2s})g(\sigma_{2s}))\,ds)^k\leq K\Delta_n^{k+1}$ for
some $k>2$ and constant $K>0$, using H\"{o}lder's inequality we have
that $C_3$ is also asymptotically negligible. This proves the
asymptotic negligibility of the term $A_4$.

We are left with proving asymptotic negligibility of $A_5$. We start
with some preliminary results that we will make use of. We have for
$0<p<\beta\wedge 1$
%e6.32 ###
\begin{equation}\label{eq:sv_a3_1}
\mathbb{E}_{i-1}^n\biggl|\Delta_n^{-1/\beta}\int_{(i-1)\Delta_n}^{i\Delta_n}\bigl(m_{ds}-m_{d,(i-1)\Delta_n}\bigr)\,ds\biggr|^p
\leq
K\Delta_n^{3p/2-p/\beta},
\end{equation}
where we made use of H\"{o}lder's inequality and the fact that $m_{ds}$
is an \Ito semimartingale with bounded coefficients and therefore
$\mathbb{E}|m_{ds}-m_{d,(i-1)\Delta_n}|\leq K|s-(i-1)\Delta_n|^{1/2}$
for $s\in[(i-1)\Delta_n,i\Delta_n)$. Similarly for $p\leq\beta$ and
arbitrary $\varepsilon>0$
%e6.33 ###
\begin{eqnarray}\label{eq:sv_a3_2}
&&
\mathbb{E}_{i-1}^n\biggl|\Delta_n^{-1/\beta}\int_{(i-1)\Delta_n}^{i\Delta_n}\int_{\mathbb{R}}(\sigma_{2s-}-\overline{\sigma}_{2s-})\kappa(x)\widetilde{\mu}(ds,dx)\nonumber
\\
&&\hphantom{\mathbb{E}_{i-1}^n\biggl|}{}+
\Delta_n^{-1/\beta}\int_{(i-1)\Delta_n}^{i\Delta_n}\int_{\mathbb{R}}(\sigma_{2s-}-\overline{\sigma}_{2s-})\kappa'(x)\mu(ds,dx)\biggr|^p
\\
&&\qquad\leq
K\Delta_n^{p/\beta-\varepsilon},\nonumber
\end{eqnarray}
where we made use of H\"{o}lder's inequality, the
Burkholder--Davis--Gundy inequality (recall $\beta>1$) and
\refeq{eq:sv_a2_6}.

Further, for some deterministic sequence $\varepsilon_n\downarrow 0$
denote
\[
S_{i}^n:=\{\omega\dvtx\Delta_n^{-1/\beta}|\Delta_i^n\overline{X}|
>
\varepsilon_n \cap \Delta_n^{-1/\beta}|\Delta_i^n X-\Delta_i^n \overline{X}|
<
0.5\varepsilon_n\}.
\]
 Then we can apply the
result in \refeq{eq:an_3b} to get for any $\alpha,\alpha'\in (0,1)$
%e6.34 ###
\begin{equation}\label{eq:sv_a3_3}
\mathbb{P}_{i-1}^n(\Delta_n^{-1/\beta}|\Delta_i^n \overline{X}|\leq\varepsilon_n)
\leq
K\biggl(\varepsilon_n^{\alpha}+\Delta_n^{(1-1/\beta)\alpha}+\frac{\Delta_n^{1-\beta'/\beta-\alpha'}}{\varepsilon_n^{\beta'}}\biggr).
\end{equation}
Similarly using the same arguments as above and \refeq{eq:an_3d}, we
get for $\varepsilon_n\downarrow 0$, some $\alpha>0$, and any $\alpha'>0$
\begin{eqnarray}\label{eq:sv_a3_4}
&&
\mathbb{E}_{i-1}^n\bigl(|\Delta_n^{-1/\beta}\Delta_i^n\overline{X}|^{-\alpha}1_{\{|\Delta_n^{-1/\beta}\Delta_i^n\overline{X}|>\varepsilon_n\}}\bigr)\nonumber
\\[-8pt]\\[-8pt]
&&\qquad\leq
K\biggl(\varepsilon_n^{(1-\alpha)\wedge 0-\alpha'}+\frac{\Delta_n^{1-\beta'/\beta-\alpha'}}{\varepsilon_n^{\alpha+\beta'}}\biggr).\nonumber
\end{eqnarray}
Finally, using \refeq{eq:sv_a3_1} and \refeq{eq:sv_a3_2}, we get for
any $\alpha'>0$
%e6.35 ###
\begin{equation}\label{eq:sv_a3_5}
\mathbb{P}_{i-1}^n(\Delta_n^{-1/\beta}|\Delta_i^n X-\Delta_i^n\overline{X}|\geq 0.5\varepsilon_n)
\leq
K\frac{\Delta_n^{1-\alpha'}}{\varepsilon_n^{\beta}}.
\end{equation}
We are now ready to prove the asymptotic negligibility of $A_5$. We can
make the following decomposition using a Taylor expansion on the set
$S_i^n$:
\begin{eqnarray*}%\label{eq:sv_a3_6}
&&
|\Delta_n^{-1/\beta}\Delta_i^n X|^p-|\Delta_n^{-1/\beta}\Delta_i^n\overline{X}|^p
\\
&&\qquad=
g(\Delta_n^{-1/\beta}\Delta_i^n X^{*})(\Delta_n^{-1/\beta}\Delta_i^n X-\Delta_n^{-1/\beta}\Delta_i^n \overline{X})1_{S_i^n}
\\
&&\qquad\quad{}+
(|\Delta_n^{-1/\beta}\Delta_i^n X|^p-|\Delta_n^{-1/\beta}\Delta_i^n \overline{X}|^p)1_{(S_i^n)^c},
\end{eqnarray*}
where $\Delta_i^n X^{*}$ is between $\Delta_i^n X$ and $\Delta_i^n
\overline{X}$ and recall $g(x)=p\mbox{ }\operatorname{sign}\{x\}|x|^{p-1}$.
Then using the definition of the set $S_i^n$ we have
$|\Delta_n^{-1/\beta} \Delta_i^n X^{*}|\geq 0.5|\Delta_n^{-1/\beta}
\Delta_i^n \overline{X}|$. Therefore, using the definition of  the
function $g(\cdot)$, it clearly suffices to show
\begin{eqnarray}
\label{eq:sv_a3_7}\qquad T_1
&:=&
\Delta_n^{-1/2}\mathbb{E}_{i-1}^n(|\Delta_n^{-1/\beta}\Delta_i^n\overline{X}|^{p-1}|\Delta_n^{-1/\beta}\Delta_i^n X-\Delta_n^{-1/\beta}\Delta_i^n
\overline{X}|1_{S_i^n})\nonumber
\\[-8pt]\\[-8pt]
&\leq&
K\Delta_n^{\alpha'},\nonumber
\\
\label{eq:sv_a3_8}T_2
&:=&
\Delta_n^{-1/2}\bigl|\mathbb{E}_{i-1}^n\bigl\{(|\Delta_n^{-1/\beta}\Delta_i^n X|^p-|\Delta_n^{-1/\beta}\Delta_i^n\overline{X}|^p)1_{(S_i^n)^c}\bigr\}\bigr|\nonumber
\\[-8pt]\\[-8pt]
&\leq&
K\Delta_n^{\alpha'}\nonumber
\end{eqnarray}
for some $\alpha'>0$. Setting $\varepsilon_n=\Delta_n^x$ for some $x>0$,
we can use the H\"{o}lder inequality to bound $T_1$
\begin{eqnarray}\label{eq:sv_a3_9_a}
T_1
&\leq&
\Delta_n^{-1/2}\bigl(\mathbb{E}_{i-1}^n|\Delta_n^{-1/\beta}\Delta_i^n\overline{X}|^{(p-1)\beta/(\beta-1)}1_{S_i^n}\bigr)^{(\beta-1)/\beta}\nonumber
\\[-8pt]\\[-8pt]
&&{}\times
(\mathbb{E}_{i-1}^n|\Delta_n^{-1/\beta}\Delta_i^nX-\Delta_n^{-1/\beta}\Delta_i^n\overline{X}|^{\beta})^{1/\beta}.\nonumber
\end{eqnarray}
Then using the bounds in \refeq{eq:sv_a3_1}, \refeq{eq:sv_a3_2} and
\refeq{eq:sv_a3_4} we get
\begin{eqnarray}\label{eq:sv_a3_9}
T_1
&\leq&
K\Delta_n^{1/\beta-1/2-\alpha'}\bigl(\Delta_n^{x(p-1+(\beta-1)/\beta)\wedge 0}\nonumber
\\[-8pt]\\[-8pt]
&&\hphantom{K\Delta_n^{1/\beta-1/2-\alpha'}\bigl(}{}+
\Delta_n^{(1-\beta'/\beta)(\beta-1)/\beta-x(1-p)-x\beta'(\beta-1)/\beta}\bigr)\nonumber
\end{eqnarray}
for some $\alpha'>0$. Similarly for $T_2$ we can use H\"{o}lder's
inequality to get
\begin{eqnarray}\label{eq:sv_a3_10_a}
T_2
&\leq&
\Delta_n^{-1/2}\bigl(\mathbb{E}_{i-1}^n\bigl||\Delta_n^{-1/\beta}\Delta_i^n X|^{p}-|\Delta_n^{-1/\beta}\Delta_i^n\overline{X}|^{p}\bigr|^{\beta/p}\bigr)^{p/\beta}\nonumber
\\[-8pt]\\[-8pt]
&&{}\times
(\mathbb{P}_{i-1}^n((S_i^n)^c))^{1-p/\beta}.\nonumber
\end{eqnarray}
Then using the bounds in \refeq{eq:sv_a3_1}, \refeq{eq:sv_a3_2},
\refeq{eq:sv_a3_3} and \refeq{eq:sv_a3_5} we get
\begin{eqnarray}\label{eq:sv_a3_10}
\qquad T_2
&\leq&
K\Delta_n^{p/\beta-1/2-\alpha'}\bigl(\Delta_n^{(1-p/\beta)x}+\Delta_n^{(1-1/\beta)(1-p/\beta)}\nonumber
\\[-8pt]\\[-8pt]
&&\hphantom{K\Delta_n^{p/\beta-1/2-\alpha'}\bigl(}{}+
\Delta_n^{(1-\beta'/\beta)(1-p/\beta)-(1-p/\beta)\beta'x}+\Delta_n^{(1-p/\beta)(1-x\beta)}\bigr)\hspace*{-8pt}\nonumber
\end{eqnarray}
for some $\alpha'>0$. Finally, we can make use of the restrictions on
$p$ and $\beta'$ to pick $x>\frac{\beta-2p}{2(\beta-p)}$ for which
\refeq{eq:sv_a3_7} and \refeq{eq:sv_a3_8} will be fulfilled.
\end{pf*}

%s6.7 ###
\subsection{\texorpdfstring{Proof of Theorem~\protect\ref{thms:sv_est}}{Proof of Theorem 4.2.}}

The proof follows directly from the fact that under the conditions of
the theorem: (1) the functions $\mu_p(\beta)$ and $\mu_{p,p}(\beta)$
are continuous both in $\beta$ and $p$; (2)
$\hat{\beta}_{X,T}^{ts}$ is consistent for $\beta$; (3)
$V_T(X,p,\Delta_n)$ converges uniformly in $p$ (after scaling
appropriately).

\section*{Acknowledgments}
We would like to thank Per Mykland, Neil Shephard and particularly Jean
Jacod for many helpful comments and suggestions. We also thank an
anonymous referee for careful reading and constructive comments on the
paper.

% imsref loaded by arune.pranskunaite, 2010-09-14 15:32:56

\printaddresses

\end{document}